\documentclass[a4paper,11pt]{article}
\usepackage{amsmath,latexsym,amssymb,fullpage}

\usepackage{chngcntr}

\newtheorem{theorem}{Theorem}[section]
\newtheorem{corollary}[theorem]{Corollary}
\newtheorem{observation}[theorem]{Observation}
\newtheorem{lemma}[theorem]{Lemma}
\newtheorem{definition}[theorem]{Definition}
\newtheorem{remark}[theorem]{Remark}
\newtheorem{proposition}[theorem]{Proposition}
\newtheorem{problem}{Problem}

\newenvironment{proof}[1][Proof:]
 {\noindent {\bf #1} }{\hspace*{0pt}\hfill$\Box$ \bigskip}

\usepackage{enumerate}
\usepackage[utf8]{inputenc}

\newif\ifskip
\skiptrue

\newcommand{\Z}{\mathbb{Z}}

\newcommand{\N}{\mathbb{N}}

\newcommand{\MM}{\mathbb{M}}
\newcommand{\cC}{\mathcal{C}}

\newcommand{\cA}{\mathcal{A}}

\newcommand{\MSOL}{\mathrm{MSOL}}
\newcommand{\CMSOL}{\mathrm{CMSOL}}
\newcommand{\FOL}{\mathrm{FOL}}

\newcommand{\C}{\mathrm{C}}

\title{Extensions and Limits of the Specker-Blatter Theorem}
\author{Eldar Fischer and Johann A. Makowsky \\ \ \\
Faculty of Computer Science\\ Israel Institute of Technology, Haifa, Israel}

\begin{document}

\maketitle

\begin{abstract}
\large
The original Specker-Blatter Theorem (1983) was formulated for classes of structures $\mathcal{C}$ 
of one or several binary relations definable in Monadic Second Order Logic MSOL. 
It states that the number of such structures on the set $[n]$ is modularly C-finite (MC-finite). 
In previous work we extended this to structures definable in CMSOL, MSOL extended with modular counting quantifiers. 
The first author also showed that the Specker-Blatter Theorem does not hold for one quaternary relation (2003).

If the vocabulary allows a constant symbol $c$, there are $n$ possible interpretations on $[n]$ for $c$.
We say that a constant $c$ is {\em hard-wired} if $c$ is always interpreted by the same element $j \in [n]$.
In this paper we show:
\begin{enumerate}[(i)]
\item 
The Specker-Blatter Theorem also holds for CMSOL when hard-wired constants are allowed. 
The proof method of Specker and Blatter does not work in this case.
\item 
The Specker-Blatter Theorem does not hold already for $\mathcal{C}$ with one ternary 
relation definable in First Order Logic FOL. This was left open since 1983.
\end{enumerate}

Using hard-wired constants allows us to show MC-finiteness of counting functions of
various restricted partition functions which were not known to be MC-finite till now. 
Among them we have the restricted Bell numbers $B_{r,A}$, 
restricted Stirling numbers of the second kind $S_{r,A}$ or restricted Lah-numbers $L_{r,A}$. 
Here $r$ is an non-negative integer and $A$ is an ultimately periodic set of non-negative integers.

\end{abstract}

\newpage
\small
\tableofcontents
\newpage
\normalsize

\part{\Large The Specker-Blatter Theorem}
\label{part-1}
\section{Introduction}

A sequence of natural numbers $s(n)$ is {\em C-finite} if it satisfies a linear recurrence relation with constant coefficients.
$s(n)$ is MC-finite if it satisfies a linear recurrence relation with constant coefficients
modulo $m$ for each $m$ separately.
A C-finite sequence $s(n)$ must have limited growth: $s(n) \leq 2^{cn}$ for some constant $c$.
No such bound exists for MC-finite sequences: for every monotone increasing sequence $s(n)$ the sequence
$s'(n)= n! s(n)$ is MC-finite.

A typical example of an C-finite sequence is the sequence $f(n)$ of Fibonacci numbers.
A typical example of an MC-finite sequence which is not C-finite is the sequence $B(n)$ of Bell numbers.
The Bell number $B(n)$ counts the number of partitions of the set $[n]$ of the numbers
$\{1, 2, \ldots, n\}$. Let $Eq(n)$ be number of equivalence relations over $[n]$. Clearly, $B(n)= Eq(n)$.
Let $Eq_2(n)$  be the number of equivalence relations on $[n]$ with exactly two equivalence classes of the same size.
$Eq_2(n)$ is not MC-finite since the value of $Eq_2(n)$ is odd iff $n$ is an even power of $2$, see \cite{pr:BlatterSpecker84}.

In \cite{pfeiffer2004counting} G. Pfeiffer discusses counting other transitive relations besides $Eq(n)$,
in particular, partial orders $PO(n)$, quasi-orders (aka preorders) $QO(n)$ and just transitive relations $Tr(n)$.
Using a growth argument one can see that none of these functions is C-finite.
It follows directly from the Specker-Blatter Theorem stated below, see Corollary \ref{cor1},
that $PO(n), QO(n), Tr(n)$ are MC-finite.
However, to the best of our knowledge, this has not been stated in the literature.
This may be due to the fact that no explicit formula for these functions are known.
The Specker-Blatter Theorem establishes MC-finiteness even in the absence of explicit formulas.
It derives MC-finiteness solely from the assumption that $\cC$ is definable in Monadic Second Order Logic
($\MSOL$), or in $\MSOL$ augmented by modular counting quantifiers ($\CMSOL$).
The present paper grew out of our study of modular recurrence relations for
restricted partition functions, \cite{ar:FischerMakowskyRakita22}.

\section{Background in logic}
We generally follow the notation of \cite{modelbook}, and assume basic knowledge of model theory. 
Standard texts for Finite Model Theory are 
\cite{modelbook,bk:Libkin2004}.
In the following, we always refer to a set $\bar R=\{R_1,\ldots,R_{\ell_{\bar R}}\}$ of distinct binary relation symbols, 
a set $\bar a=\{a_1,\ldots,a_{\ell_{\bar a}}\}$ of distinct constant symbols, and so on. 
By $a\in \bar a$ we mean that there exists $1\leq i\leq\ell_{\bar a}$ for which $a=a_i$. 
We also use the shorthand $[n]=\{1,\ldots,n\}$.

Let $\tau = \bar{R} \cup \bar{a}$ be a vocabulary, i.e., a set of non-logical constants.
We denote by $\FOL(\tau)$ the set of first order formulas with its non-logical constants in $\tau$.
If $\tau$ is clear from the context, we omit it.
We denote by $\MSOL(\tau)$ the set of Monadic Second Order Logic, obtained from $\FOL$ by allowing
unary relation variables and quantification over them.
The logic $\CMSOL$ is obtained from $\MSOL$ by allowing also quantification of the form $\C_{m,a}x \phi(x)$,
which are interpreted by 
$$
\mathcal{A} \models \C_{m,a}x \phi(x) \text{   iff   } |\{ a \in A : \phi(a) \}| \equiv a \mod{m}
$$

\section{The original Specker-Blatter Theorem}
Let $\phi_E$ be the formula in First Order Logic ($\FOL$) which says that $E(x,y)$ is an equivalence relation.
$Eq(n)$ can be written as
$$
Eq(n) = |\{ E \subseteq [n]^2 : ([n], E) \models \phi_E \}|.
$$
$PO(n), QO(n), Tr(n)$ can be written in a similar way.

The original Specker-Blatter Theorem from 1981, \cite{blatter1981nombre,SB,pr:BlatterSpecker84,specker1990application}, 
gives a general criterion for certain integer sequences
to be MC-finite. 
Let $\bar R = \{R_1, \ldots, R_{\ell_{\bar R}} \}$ be a finite set of relation symbols
of arities $\rho_1,\ldots,\rho_{\ell_{\bar R}}$ respectively, and $\phi$ be a formula of Monadic Second Order Logic ($\MSOL$)
using relation symbols from $R$ without free variables.

Let $Sp_{\phi}(n)$ be the number of ways we can interpret the relation symbols in $R$ on $[n]$ 
such that the resulting structures where $A_i$ is the interpretation of $R_i$ satisfies $\phi$.
Formally
$$
Sp_{\phi}(n)=
|\{
A_i \subseteq [n]^{\rho_i}, i \leq m : ([n], A_1. \ldots, A_m) \models \phi
\}|
$$

\begin{theorem}
Let $\bar{R}$ be a finite set of binary relations and $\phi$ be 
a formula of $\MSOL(\bar{R})$ using relation symbols in $\bar{R}$.
Then the sequence  $Sp_{\phi}(n)$ is MC-finite.
\end{theorem}

\begin{corollary}
\label{cor1}
The sequences counting the number of partial orders $PO(n)$, quasi-orders $QO(n$), and
transitive relations $Tr(n)$ on $[n]$, are MC-finite.
\end{corollary}
Some recent applications of the Specker-Blatter Theorem can be found in \cite{kotek2021tutte,ar:FischerMakowskyRakita22}.

\section{Previous limitations and extensions}
Limitations and extensions of the Specker-Blatter Theorem have been previously discussed in
\cite{pr:FischerMakowsky03,sbplus}.

It is well known that Eulerian graphs and regular graphs of even degree are not definable in $\MSOL$,
but they are definable in $\CMSOL$.
In
\cite{pr:FischerMakowsky03},
the Specker-Blatter Theorem was shown to hold also for $\CMSOL$.

It follows in particular that $Eul(n)$, which counts the number of Eulerian graphs over $[n]$ (i.e.\ connected graphs all of whose degrees are even), is also MC-finite.

In 
\cite{noquad}
the first author showed in 2002 that the Specker-Blatter Theorem does not hold for
quaternary relations.

\begin{theorem}[E.\ Fischer, 2002]
\label{th:quaternary}
There is an $\FOL$-sentence with only one quaternary relation symbol $\phi_1$, such that
$Sp_{\phi_1}(n)$ is not an MC-sequence.
\end{theorem}

The question of whether Specker-Blatter Theorem holds in the presence of a ternary relation symbols
remained open. 

\section{Substitution rank}

The proof of the Specker-Blatter Theorem has two parts, both of which use the assertion that $\tau= \bar{R}$ contains
only binary relation symbols. One is purely combinatorial and uses an
operation between $\bar{R}$-structures called {\em substitution}.
A pointed $\bar{R}$-structure is an $\bar{R}$-structure $\mathcal{A}_a =([n_1], A_1, \ldots, A_{\mu}, a)$ with an additional
distinguished point $a \in [n]$. Given a pointed $\bar{R}$-structure ${\mathcal{A}_a}_1$ and 
an $\bar{R}$-structure $\mathcal{A}_2$ with universe $[n_2]$ we define $\cA =Subst(\cA_1,a,\cA_2)$ as follows:
\begin{enumerate}[(i)]
\item
The universe $A$ of $\cA$ is the disjoint union of $A_1$ and $A_2$ with the point $a$ removed.
It can be assumed to be the set $[n_1+n_2-1]$.
\item
The binary relations 
are defined such that $\cA_2$ is a module in $\cA$, i.e., for $u \in A_1-\{a\}$ and $v \in A_2$
and $R \in \bar{R}$, the relation $R(u,v)$ holds in $\cA =Subst(\cA_1,a,\cA_2)$ iff $R(u,a)$ holds in
${\mathcal{A}_a}_1$.
\end{enumerate}


Now let $(\cA_a)_i$ 
be an enumeration of all pointed $R$-structures
and $(\cA)_j$
be an enumeration of all (non-pointed) $R$-structures.
Let $\cC$ be a class of finite of $R$-structures.
We define the infinite $0/1$-matrix $\MM(P)$ for $P$ as follows.
\begin{enumerate}[(i)]
\item
The rows are labeled by the pointed structures  $(\cA_a)_i$.
\item
The columns are labeled by the non-pointed structures  $(\cA)_j$.
\item
$\MM(\cC)_{i,j} =1$ iff $Subst((\cA_a)_i,a,\cA_j) \in \cC$.
\end{enumerate}
The substitution rank of $\cC$ is the rank of the matrix $\MM(\cC)$ over the ring $\Z_{p^q}$ under discussion (for working with $\Z_m$ for any integer $m$, we work with each of its maximum prime power divisors separately).
We define
$$
Sp_{\cC}(n)=
|\{
A_i \subseteq [n]^{\rho_i}, i \leq m : ([n], A_1. \ldots, A_m) \in \cC
\}|
$$
\begin{theorem}
\label{th:rank}
Let $\bar{R}$ be a finite set of binary relations and 
$\cC$ be a class of finite $\bar{R}$-structures whose substitution rank is finite under any $\Z_{p_q}$.
Then the sequence  $Sp_{\cC}(n)$ is MC-finite.
\end{theorem}
This version of the Specker-Blatter Theorem applies to uncountable many properties
of $\bar{R}$-structures, whereas there are only countably many such properties definable in $\MSOL$.

\begin{theorem}
\label{th:msol}
Let $\bar{R}$ be a finite set of binary relations and 
$\cC$ be 
a finite class of $R$-structures defined by an
$R$-sentences $\phi$ in $\MSOL$.
Then the substitution rank of $\cC$ is finite.
\end{theorem}

In
\cite{pr:FischerMakowsky03} it is shown that this still holds if $\MSOL$ is replaced by $\CMSOL$.

If the vocabulary allows a constant symbol $c$, there are $n$ possible interpretations on $[n]$ for $c$.
We say that a constant $c$ is {\em hard-wired} if $c$ is always interpreted by the same element $j \in [n]$.
\begin{remark}
The definition of $\cA =Subst(\cA_1,a,\cA_2)$ relies on the fact that all relation symbols are binary.
In particular it is not clear how to handle hard-wired constants.
A unary relation $U(x)$ can be replaced by a binary relation $R(x,x)$ with
$\forall x (U(x) \leftrightarrow R(x,x))$ and $\forall x,y ((x \neq y) \rightarrow \neg R(x,y))$.
For relation symbols of arity $\geq 3$ there are various options for defining the substitution 
$\cA =Subst(\cA_1,a,\cA_2)$, but none seems to work to prove Theorem \ref{th:rank}.
\end{remark}

\section{Main new results}
\label{newresults}
The Bell numbers $B(n)$ and the Stirling numbers of the second kind $S_k(n)$ for fixed $k$
can be shown to be MC-finite using the Specker-Blatter Theorem.
A. Broder in 1984, \cite{broder1984r}, introduced the restricted Bell numbers $B_r(n)$ and the
restricted Stirling numbers of the second kind $S_{k,r}(n)$.
Let $r \in \N^+$. $S_{k,r}(n)$ is defined as the number of set partitions of $[r+n]$ into $k+r$ blocks with
the additional condition that the first $r$ elements are in distinct blocks. $B_r(n)$ is defined as
$$
B_r(n) = \sum_k S_{k,r}(n)
$$
The class of equivalence relations on $[r+n]$ where the first $r$ elements are in different equivalence classes
is definable in $\FOL$ with one binary relation and $r$ hard-wired constants.
The Specker-Blatter Theorem does directly apply to this case. In \cite{ar:FischerMakowskyRakita22}
it is shown how to circumvent this obstacle in the case of one equivalence relation.
It followed that both $S_{k,r}(n)$ and $B_r(n)$ are MC-finite.

In this paper we prove a more general theorem:

\begin{theorem}[Elimination of hard-wired constants]
\label{th:main-1} ~~~
\begin{enumerate}[(i)]
\item
Let $\tau$ consist of a  finite set of constant symbols $\bar{a}$, unary relations symbols $\bar{U}$, and binary relation symbols $\bar{R}$.
For every class $\cC$ of $\tau$-structures there exist classes $\cC_1, \ldots, \cC_r$
of $\tau'$-structures, where $\tau'$-contains only a finite number $r(\bar{a},\bar{U},\bar{R})$ of binary relation symbols, such that
$$
Sp_{\cC}(n) = \sum_{i=1}^r Sp_{\cC_i}(n).
$$
Equality here is not modular.
\item
Furthermore, of $\cC$ is $\FOL$-definable ($\MSOL$-definable, $\CMSOL$-definable), so are the $\cC_i$.
\end{enumerate}
\end{theorem}

\begin{corollary}\label{cr:spc}
Let $\tau$ consist of a  finite set of constant symbols $\bar{a}$, unary relations symbols $\bar{U}$, and
binary relation symbols $\bar{R}$, let $\cC$ be class of finite $\tau$-structures definable
in $\CMSOL$. Then the sequence $Sp_{\cC}(n)$ is  MC-finite.
\end{corollary}
The proof is given in Part \ref{part-2}

We have seen in Theorem \ref{th:quaternary}
that the Specker-Blatter Theorem does not hold for a 
single quaternary relation. 
The question of whether Specker-Blatter Theorem holds in the presence of a single ternary relation symbol
remained open. Our second main result here answers this question. The proof is given in Part \ref{part-3}

\begin{theorem}[Ternary Counter-Example]
\label{th:ternary}
There is a $\FOL$-sentence $\phi$ with only one ternary relation symbol,
such that $Sp_{\phi}(n)$ is not an MC-sequence.
\end{theorem}

The proof of Theorem \ref{th:ternary} first produces a sentence $\psi$ which also uses one symbol for a hard-wired constant.
This will be shown in Section \ref{se:ternary}.
To construct $\phi$ we need a way to eliminate hard-wired constant symbols.
For this we use the technique developed in Part \ref{part-2}, which converts $\psi$ to a sentence with a single ternary relation and a constant number of lower arity relations. In Section \ref{s:s} we show how for this particular sentence we can then get rid of the added relations, leaving us with only the ternary relation.

We conclude the paper with Part \ref{part-3}, Section \ref{se:conclu}, containing a summary and open problems.


\part{\Large The ephemeral role of hard-coded constants}
\label{part-2}
\section{Introduction}


In the following we consider extending the language with ``hard-coded'' constants. Specifically, assume that we have a class $\mathcal C$ that is defined by a sentence $\phi$ involving a set of constant symbols $\bar a$, unary symbols $\bar U$ and binary symbols $\bar R$. The function $f_{\mathcal C}(n)$ is defined as the number of models over the universe $[n+\ell_{\bar a}]$ which satisfy $\phi$, for which $a_i$ is interpreted as $n+i$ for all $i\in [\ell_{\bar a}]$. Note the distinction from the non-hard-coded setting, where we would have had to also count the possible interpretations of the constants.

Our main result is an expression for the function $f_{\mathcal C}(n)$ (when constants are allowed) that is based on counting functions for classes that do not utilize constants. We first show this reduction for languages using only unary and binary relations, and then show how to extend it for higher arities. The reduction preserves many of the common logics, in particular an FOL expression would be reduced to functions involving FOL expressions, and so on. On the positive side, this extends the Specker-Blatter theorem to languages involving hard-coded constants, allowing modular ultimate periodicity proofs of new functions.

On the negative side, we construct an FOL statement involving a ternary relation and a constant that (building on the quaternary example) produces a function that is not ultimately periodic modulo $2$. By the reduction, this means that there exists such a statement not involving a constant. We later show how ``trim down'' this statement so it uses a single ternary relation and no other (lower arity) relations, completely closing the knowledge gap.

\section{Proving the reduction}\label{s:r}

In this section we prove Theorem \ref{th:main-1}.
For convenience we state it again as Theorem \ref{th:r}.

\begin{theorem}[Reducing model counts to the case without constants]\label{th:r}
For any class $\mathcal C$ defined by an FOL (resp.\ MSOL, CMSOL) sentence involving a set of constant symbols $\bar a$, unary symbols $\bar U$ and binary symbols $\bar R$, there exist classes $\mathcal C_1,\ldots,\mathcal C_r$ (where $r$ depends on the original language), definable by FOL (resp.\ MSOL, CMSOL) sentences involving $\bar U'$ (which contains $\bar U$), $\bar R$ and no constants, satisfying $f_{\mathcal C}(n)=\sum_{i=1}^rf_{\mathcal C_i}(n)$ for all $n\in\mathbb N$.
\end{theorem}

Later on, following similar lines, we streamline this theorem to use a single target class (a many-one reduction), and then extend to classes involving higher rank relations. The following is the immediate corollary it produces for the Specker-Blatter Theorem.

\begin{corollary}[Extended Specker-Blatter Theorem]\label{cor:sb}
For a class $\mathcal C$ definable in CMSOL with (hard-coded) constants, unary and binary relation symbols only, 
the function $f_{\mathcal C}$ is MC-finite.
\end{corollary}

Theorem \ref{th:r} is proved by induction over the number of constants. The basis, $\ell_{\bar a}=0$, is trivial (with $\bar U'=\bar U$, $r=1$ and $\mathcal C_1=\mathcal C$). The induction step is provided by the following.

\begin{lemma}[Removing a single constant]\label{lm:scr}
For any class $\mathcal C$ defined by an FOL (resp.\ MSOL, CMSOL) sentence involving a set of constant symbols $\bar a$ with $\ell_{\bar a}>0$, unary symbols $\bar U$ and binary symbols $\bar R$, there exist classes $\mathcal C_1,\ldots,\mathcal C_r$ (where $r$ depends on the original language), definable by FOL (resp.\ MSOL, CMSOL) sentences over the language $(\bar a',\bar U',\bar R')$, where $\bar a'=\bar a\setminus\{a_{\ell_{\bar a}}\}$, $\bar U'=\bar U\cup\bar I\cup\bar O$ where $\ell_{\bar I}=\ell_{\bar O}=\ell_{\bar R}$, and $\bar R'=\bar R$, satisfying $f_{\mathcal C}(n)=\sum_{i=1}^rf_{\mathcal C_i}(n)$ for all $n\in\mathbb N$.
\end{lemma}

The main idea in the proof of this lemma is to encode the ``interaction'' of the constant $a_{\ell_{\bar a}}$ with the rest of the universe using the additional unary relations. For every $i\in [\ell_{\bar R}]$, we will use the new relation $I_i$ to hold every $x\neq a_{\ell_{\bar a}}$ for which $(x,a)$ was in $R_i$, and the relation $O_i$ to hold every $x\neq a_{\ell_{\bar a}}$ for which $(a,x)$ was in $R_i$.

We cannot directly keep track whether $(a,a)$ was in $R_i$, or whether $a$ was in $U_i$ for $i\in [\ell_{\bar U}]$, so we count the number of models for each of these options separately. This sets $r=2^{\ell_{\bar U}+\ell_{\bar R}}$. Instead of a running index, we index each such option with a set $\mathfrak U\subseteq [\ell_{\bar U}]$ denoting which of the relations in $\bar U$ include the constant to be removed $a=a_{\ell_{\bar a}}$, and a set $\mathfrak R\subseteq [\ell_{\bar R}]$ denoting which of the relations in $\bar R$ include $(a,a)$. Using these we can define the case where a model $\mathfrak N$ over the language $(\bar a',\bar U',\bar R)$ with universe $[n+\ell_{\bar a}-1]$ corresponds (along with $\mathfrak U$ and $\mathfrak R$) to an ``original model'' $\mathfrak M$ with universe $[n+\ell_{\bar a}]$ over the original language.

\begin{definition}\label{df:corr}
Given a model $\mathfrak M$ over the language $(\bar a,\bar U,\bar R)$ with universe $[n+\ell_{\bar a}]$, a model $\mathfrak N$ over the language $(\bar a',\bar U',\bar R)$ with universe $[n+\ell_{\bar a}-1]$, and sets $\mathfrak U\subseteq [\ell_{\bar U}]$ and $\mathfrak R\subseteq [\ell_{\bar R}]$, where (as always) in both models every constant $a_i$ is interpreted to be $n+i$, we say that $(\mathfrak N,\mathfrak U,\mathfrak R)$ {\em correspond} to $\mathfrak M$ if the following hold.
\begin{itemize}
\item For every $U\in \bar U$ and $x\in [n+\ell_{\bar a}-1]$, we have $\mathfrak N\models U(x)$ if and only if $\mathfrak M\models U(x)$.
\item For every $i\in [\ell_{\bar U}]$, we have $i\in\mathfrak U$ if and only if $\mathfrak M\models U_i(a)$.
\item For every $R\in \bar R$ and $x,y\in [n+\ell_{\bar a}-1]$, we have $\mathfrak N\models R(x,y)$ if and only if $\mathfrak M\models R(x,y)$.
\item For every $i\in [\ell_{\bar R}]$ and $x\in [n+\ell_{\bar a}-1]$, we have $\mathfrak N\models I_i(x)$ if and only if $\mathfrak M\models R_i(x,a)$.
\item For every $i\in [\ell_{\bar R}]$ and $x\in [n+\ell_{\bar a}-1]$, we have $\mathfrak N\models O_i(x)$ if and only if $\mathfrak M\models R_i(a,x)$.
\item For every $i\in [\ell_{\bar R}]$, we have $i\in\mathfrak R$ if and only if $\mathfrak M\models R_i(a,a)$.
\end{itemize}
\end{definition}

It is important to note, for the purpose of counting, the following observation.

\begin{observation}\label{ob:bij}
Definition \ref{df:corr} provides a bijection between the set of possible models $\mathfrak M$ over the universe $[n+\ell_{\bar a}]$ (where the constants are interpreted as in Definition \ref{df:corr}), and the set of possible triples $(\mathfrak N,\mathfrak U,\mathfrak R)$ where $\mathfrak N$ is a model over $[n+\ell_{\bar a}-1]$  (where the constants are interpreted as in Definition \ref{df:corr}) and $\mathfrak U\subseteq [\ell_{\bar U}]$ and $\mathfrak R\subseteq [\ell_{\bar R}]$.
\end{observation}

Suppose we are given an expression $\phi(\bar x)$ where $\bar x=\{x_1,\ldots,x_{\ell_{\bar x}}\}$ is a set of variable symbols over the language $(\bar a,\bar U,\bar R)$, as well as a set $\mathfrak U\subseteq [\ell_{\bar U}]$ and a set $\mathfrak R\subseteq [\ell_{\bar R}]$. We will construct, by induction over the structure of $\phi$, several expressions, where one of which is an expression $\phi'_{\mathfrak U,\mathfrak R}(\bar x)$ over the language $(\bar a',\bar U',\bar R)$. It will be constructed so that for any $\mathfrak M$ over the language $(\bar a,\bar U,\bar R)$ with universe $[n+\ell_{\bar a}]$ and $\mathfrak N$ over the language $(\bar a',\bar U',\bar R)$ with universe $[n+\ell_{\bar a}-1]$, where $(\mathfrak N,\mathfrak U,\mathfrak R)$ correspond to $\mathfrak M$, and any fixing of $x_1,\ldots,x_{\ell_{\bar x}}\in [n+\ell_{\bar a}-1]$, we will have $\mathfrak M\models\phi(\bar x)$ if and only if $\mathfrak N\models\phi'_{\mathfrak U,\mathfrak R}(\bar x)$.

Lemma \ref{lm:scr} then immediately follows from the case $\ell_{\bar x}=0$ (i.e.\ where $\phi$ is a sentence). To be precise, for a class $\mathcal C$ defined by a sentence $\phi$ over the language $(\bar a,\bar U,\bar R)$, we obtain $f_{\mathcal C}(n)=\sum_{\mathfrak U\subseteq [\ell_{\bar U}],\mathfrak R\subseteq [\ell_{\bar R}]}f_{C_{\mathfrak U,\mathfrak R}}(n)$, where $\mathcal C_{\mathfrak U,\mathfrak R}$ is the class respectively defined by $\phi'_{\mathfrak U,\mathfrak R}(\bar x)$ over the language  $(\bar a',\bar U',\bar R)$.

To sustain the induction, the above will not be enough. This is because we need to account under the model $\mathfrak N$ also for the case where some variables are ``assigned the value $a=a_{\ell_{\bar a}}$'', a value which does not exist in its universe (it exists only in that of $\mathfrak M$). We henceforth consider also a set $\mathfrak X\subseteq [\ell_{\bar x}]$, and denote the set of variable symbols $x_{\mathfrak X}=\{x_i:i\in\mathfrak X\}$. In our induction we will construct the expressions $\phi'_{\mathfrak X,\mathfrak U,\mathfrak R}(\bar x\setminus x_{\mathfrak X})$, where $\phi'_{\mathfrak U,\mathfrak R}(\bar x)$ is just the special case $\phi'_{\emptyset,\mathfrak U,\mathfrak R}(\bar x)$. With models $\mathfrak M$ and $\mathfrak N$ as above and a fixing of the variables in $\bar x\setminus x_{\mathfrak X}$, denote by $\bar x_{\mathfrak X\to a}$ the completion of this fixing to all of $\bar x$ that is obtained by fixing $x_i$ to be equal to $a$ for all $i\in\mathfrak X$. We will then have $\mathfrak M\models\phi(\bar x_{\mathfrak X\to a})$ if and only if $\mathfrak N\models\phi'_{\mathfrak X,\mathfrak U,\mathfrak R}(\bar x\setminus x_{\mathfrak X})$.

The rest of this section is concerned with the recursive definition of $\phi'_{\mathfrak X,\mathfrak U,\mathfrak R}(\bar x\setminus x_{\mathfrak X})$. There is a subsection for the base cases, a subsection for Boolean connectives, and a subsection for each class of quantifiers (first order quantifiers, counting quantifiers, and monadic second order quantifiers). In every construction we argue (at times trivially) that we keep the correspondence invariant, namely that $\mathfrak M\models\phi(\bar x_{\mathfrak X\to a})$ if and only if $\mathfrak N\models\phi'_{\mathfrak X,\mathfrak U,\mathfrak R}(\bar x\setminus x_{\mathfrak X})$ whenever $\mathfrak M$ and $(\mathfrak N,\mathfrak U,\mathfrak R)$ satisfy the correspondence condition of Definition \ref{df:corr}.

\subsection{The base constructions}\label{ss:rbase}

We use the Boolean ``true'' and ``false'' statements in the following, so for formality's sake they are also considered as atomic statements here. Clearly, if $\phi(\bar x)$ is simply the ``true'' statement $\top$ (respectively the ``false'' statement $\bot$), then setting $\phi'_{\mathfrak X,\mathfrak U,\mathfrak R}(\bar x\setminus x_{\mathfrak X})$ to $\top$ (respectively $\bot$) gives us the equivalent statement satisfying the correspondence invariant.

For $i\in [\ell_{\bar U}]$ and $j\in [\ell_{\bar x}]$, let us now consider the expression $\phi(\bar x)=U_i(x_j)$. To produce $\phi'_{\mathfrak X,\mathfrak U,\mathfrak R}(\bar x\setminus x_{\mathfrak X})$, we partition to cases according to whether $j\in\mathfrak X$. In the case where $j\notin\mathfrak X$, we simply set $\phi'_{\mathfrak X,\mathfrak U,\mathfrak R}(\bar x\setminus x_{\mathfrak X})=U_i(x_j)$ as well, which clearly satisfies the invariant for $(\mathfrak N,\mathfrak U,\mathfrak R)$ correlated with $\mathfrak M$ (recall that the ``if and only if'' condition in this case should hold when the value of $x_i$ is in $[n+\ell_{\bar a}-1]$).

Similarly, for $i\in [\ell_{\bar U}]$ and $j\in [\ell_{\bar a}-1]$, for the expression $\phi(\bar x)=U_i(a_j)$, we produce $\phi'_{\mathfrak X,\mathfrak U,\mathfrak R}(\bar x\setminus x_{\mathfrak X})=U_i(a_j)$, noting that in our setting the value of $a_j$ is guaranteed to be equal to $n+j\in [n+\ell_{\bar a}-1]$.

Now consider the expression $\phi(\bar x)=U_i(x_j)$ for the case where $x_j\in\mathfrak X$. Recall that in this case $\phi'_{\mathfrak X,\mathfrak U,\mathfrak R}(\bar x\setminus x_{\mathfrak X})$ should not depend on $x_j$. Moreover, to preserve the invariant for corresponding sets and models, $\phi'_{\mathfrak X,\mathfrak U,\mathfrak R}(\bar x\setminus x_{\mathfrak X})$ should hold if and only if $U_j(a)$ holds (recall that we use the shorthand $a=a_{\ell_{\bar a}}$ throughout). We hence define $\phi'_{\mathfrak X,\mathfrak U,\mathfrak R}(\bar x\setminus x_{\mathfrak X})$ to be $\top$ (``true'') if $i\in\mathfrak U$, and define it to be $\bot$ (``false'') if $i\notin\mathfrak U$.

The remaining case for a unary relation is the expression $\phi(\bar x)=U_i(a)$. Again, we define $\phi'_{\mathfrak X,\mathfrak U,\mathfrak R}(\bar x\setminus x_{\mathfrak X})$ to be $\top$ if $i\in\mathfrak U$, and define it to be $\bot$ if $i\notin\mathfrak U$.

We now move on to handle the atomic expressions involving a binary relation $R_i$ where $i\in [\ell_{\bar R}]$. The first case here is the one analogous to the first case we discussed involving a unary relation. Namely, it is the case where $\phi(\bar x)=R_i(x_j,x_k)$ where both $j\notin\mathfrak X$ and $k\notin\mathfrak X$. In this case we set $\phi'_{\mathfrak X,\mathfrak U,\mathfrak R}(\bar x\setminus x_{\mathfrak X})=R_i(x_j,x_k)$, and argue the same argument as above about satisfying the correspondence invariant.

The next four cases we discuss resemble the last two cases we discussed about a unary relation. Namely, these are the cases where $\phi(\bar x)=R_i(x_j,x_k)$ with $j,k\in\mathfrak X$, $\phi(\bar x)=R_i(x_j,a)$ or $\phi(\bar x)=R_i(a,x_j)$ with $j\in\mathfrak X$, and $\phi(\bar x)=R_i(a,a)$. In all these cases the resulting expression should reflect on whether $\mathfrak M\models R_i(a,a)$, which for the corresponding $(\mathfrak N,\mathfrak U,\mathfrak R)$ is handled by the set $\mathfrak R$. We hence set $\phi'_{\mathfrak X,\mathfrak U,\mathfrak R}(\bar x\setminus x_{\mathfrak X})=\top$ if $i\in\mathfrak R$, and set $\phi'_{\mathfrak X,\mathfrak U,\mathfrak R}(\bar x\setminus x_{\mathfrak X})=\bot$ if $i\notin\mathfrak R$.

Next we handle the cases where $\phi(\bar x)=R_i(x_j,x_k)$ with $j\notin\mathfrak X$ and $k\in\mathfrak X$, and $\phi(\bar x)=R_i(x_j,a)$ with $j\notin\mathfrak X$. For both this cases, for the correspondence invariant to follow we need to look at whether $\mathfrak M\models R_i(x_j,a)$, where the value of $x_j$ is in $[n+\ell_{\bar a}-1]$. For the corresponding $(\mathfrak N,\mathfrak U,\mathfrak R)$ this occurs if and only if $\mathfrak N\models I_i(x_j)$, where we recall that $I_i$ is a relation from $\bar U'\setminus\bar U$. We therefor set $\phi'_{\mathfrak X,\mathfrak U,\mathfrak R}(\bar x\setminus x_{\mathfrak X})=I_i(x_j)$ in these cases. Similarly, for the cases $\phi(\bar x)=R_i(a_j,x_k)$ and $\phi(\bar x)=R_i(a_j,a)$, where $j\in[\ell_{\bar a}-1]$ and $k\in\mathfrak X$, we set $\phi'_{\mathfrak X,\mathfrak U,\mathfrak R}(\bar x\setminus x_{\mathfrak X})=I_i(a_j)$.

Moving on to the remaining cases for a binary relation, we consider $\phi(\bar x)=R_i(x_k,x_j)$ with $j\notin\mathfrak X$ and $k\in\mathfrak X$, and $\phi(\bar x)=R_i(a,x_j)$ with $j\notin\mathfrak X$. These are analogous to the cases handled in the last paragraph, only here we use $O_i$ instead of $I_i$. We set $\phi'_{\mathfrak X,\mathfrak U,\mathfrak R}(\bar x\setminus x_{\mathfrak X})=O_i(x_j)$ in these two cases. Finally, for the cases $\phi(\bar x)=R_i(x_k,a_j)$ and $\phi(\bar x)=R_i(a,a_j)$, where $j\in[\ell_{\bar a}-1]$ and $k\in\mathfrak X$, we set $\phi'_{\mathfrak X,\mathfrak U,\mathfrak R}(\bar x\setminus x_{\mathfrak X})=O_i(a_j)$.

The final atomic formula to consider is the ``builtin relation'' of equality. We skip all cases involving only constants (e.g.\ $a_i=a_j$), since these are equivalent to $\top$ or $\bot$. We also skip cases that are equivalent by the symmetry of the equality relation to those that we discuss.

First, if $\phi(\bar x)$ is $x_i=x_j$ or $x_i=a_k$ for $i,j\notin\mathfrak X$ and $k\in [\ell_{\bar a}-1]$, then since we are dealing with values that are guaranteed to be in $[n+\ell_{\bar a}-1]$ (the universe of $\mathfrak N$), we set $\phi'_{\mathfrak X,\mathfrak U,\mathfrak R}(\bar x\setminus x_{\mathfrak X})$ respectively to $x_i=x_j$ or $x_i=a_k$ as well (so it is ``unaltered'' from $\phi(\bar x)$).

On the other hand, if $\phi(\bar x)$ is $x_i=x_j$ or $x_i=a$ for $i,j\in\mathfrak X$, then for the correspondence principle to hold, we need $\mathfrak N\models\phi'_{\mathfrak X,\mathfrak U,\mathfrak R}(\bar x\setminus x_{\mathfrak X})$ to hold if $\mathfrak M\models (a=a)$. In other words, we have to set $\phi'_{\mathfrak X,\mathfrak U,\mathfrak R}(\bar x\setminus x_{\mathfrak X})=\top$ here.

The final cases are those where $\phi(\bar x)$ is $x_i=x_j$ or $x_i=a$ for $i\notin\mathfrak X$ and $j\in\mathfrak X$. For the correspondence principle to hold, we need $\mathfrak N\models\phi'_{\mathfrak X,\mathfrak U,\mathfrak R}(\bar x\setminus x_{\mathfrak X})$ to hold if and only if $\mathfrak M\models (x_i=a)$. However, we make here the subtle yet important observation that this should occur for any value that $x_i$ can take from the universe of $\mathfrak N$, which does not include $a$. Therefor, we can (and should) set $\phi'_{\mathfrak X,\mathfrak U,\mathfrak R}(\bar x\setminus x_{\mathfrak X})=\bot$ in these cases.

\subsection{Boolean connectives}\label{ss:rbool}

Handling Boolean connectives is the most straightforward part of this construction. For example, suppose that we have $\phi(\bar x)=\neg\psi(\bar x)$ for some expression $\psi(\bar x)$, for which we have already established (by the induction hypothesis) that $\mathfrak M\models\psi(\bar x_{\mathfrak X\to a})$ if and only if $\mathfrak N\models\psi'_{\mathfrak X,\mathfrak U,\mathfrak R}(\bar x\setminus x_{\mathfrak X})$ whenever $\mathfrak M$ and $(\mathfrak N,\mathfrak U,\mathfrak R)$ correspond. Here we can clearly set $\phi'_{\mathfrak X,\mathfrak U,\mathfrak R}(\bar x\setminus x_{\mathfrak X})=\neg\psi'_{\mathfrak X,\mathfrak U,\mathfrak R}(\bar x\setminus x_{\mathfrak X})$, and obtain that $\mathfrak M\models\phi(\bar x_{\mathfrak X\to a})$ if and only if $\mathfrak N\models\phi'_{\mathfrak X,\mathfrak U,\mathfrak R}(\bar x\setminus x_{\mathfrak X})$ whenever $\mathfrak M$ and $(\mathfrak N,\mathfrak U,\mathfrak R)$ correspond.

The same idea and analysis follow for all other Boolean connectives. For example, for the expression $\phi(\bar x)=\psi_1(\bar x)\wedge\psi_2(\bar x)$, we set $\phi'_{\mathfrak X,\mathfrak U,\mathfrak R}(\bar x\setminus x_{\mathfrak X})=\psi'_{1,\mathfrak X,\mathfrak U,\mathfrak R}(\bar x\setminus x_{\mathfrak X})\wedge\psi'_{2,\mathfrak X,\mathfrak U,\mathfrak R}(\bar x\setminus x_{\mathfrak X})$.

\subsection{First order quantifiers}\label{ss:rfo}

To deal with quantifiers over variables, we make some assumptions on the structure of our expressions that can easily be justified by the appropriate variable substitutions. Namely, we require that every quantified variable is quantified only once in the expression, and is not used at all outside the scope of the quantification. In particular, this means that the set $\mathfrak X$ that appears in the subscript of our expression cannot contain a reference to the quantified variable.

For notational convenience, when $\phi(\bar x)$ is our formula, we denote by $x=x_{\ell_{\bar x}+1}$ the quantified variable. So the two cases that we consider in this subsection are the existential quantification $\phi(\bar x)=\exists_x\psi(\bar x\cup\{x\})$ and the universal quantification $\phi(\bar x)=\forall_x\psi(\bar x\cup\{x\})$, and for both of them we would like to construct a corresponding $\phi'_{\mathfrak X,\mathfrak U,\mathfrak R}(\bar x\setminus x_{\mathfrak X})$ where $\mathfrak X\subseteq [\ell_{\bar x}]$.

In the existential case, we want $\mathfrak N\models\phi'_{\mathfrak X,\mathfrak U,\mathfrak R}(\bar x\setminus x_{\mathfrak X})$ to occur whenever there is at least one value of $x$ for which $\mathfrak M\models\psi(\bar x\cup\{x\})$. For the values of $x$ within $[n+\ell_{\bar a}-1]$, by the induction hypothesis, this is covered by the expression $\exists_x\psi'_{\mathfrak X,\mathfrak U,\mathfrak R}(\bar x\cup\{x\}\setminus x_{\mathfrak X})$. However, there is one possible value of $x$ not covered in this way, and that is the value $n+\ell_{\bar a}$, which we identify with the constant $a$. But by the induction hypothesis, $\mathfrak M\models\psi(\bar x\cup\{x\})$ for $x=a$ if and only if $\mathfrak N\models\psi'_{\mathfrak X\cup\{\ell_{\bar x}+1\},\mathfrak U,\mathfrak R}(\bar x\setminus x_{\mathfrak X})$. The combined expression that satisfies the correspondence invariant is hence $\phi'_{\mathfrak X,\mathfrak U,\mathfrak R}(\bar x\setminus x_{\mathfrak X})=\exists_x\psi'_{\mathfrak X,\mathfrak U,\mathfrak R}(\bar x\cup\{x\}\setminus x_{\mathfrak X})\vee\psi'_{\mathfrak X\cup\{\ell_{\bar x}+1\},\mathfrak U,\mathfrak R}(\bar x\setminus x_{\mathfrak X})$.

The universal case follows an analogous argument, only here $\mathfrak M\models\psi(\bar x\cup\{x\})$ needs to hold for all values of $x$, those in $[n+\ell_{\bar a}-1]$ as well as the value of $a$. The combined expression is $\phi'_{\mathfrak X,\mathfrak U,\mathfrak R}(\bar x\setminus x_{\mathfrak X})=\forall_x\psi'_{\mathfrak X,\mathfrak U,\mathfrak R}(\bar x\cup\{x\}\setminus x_{\mathfrak X})\wedge\psi'_{\mathfrak X\cup\{\ell_{\bar x}+1\},\mathfrak U,\mathfrak R}(\bar x\setminus x_{\mathfrak X})$.

\subsection{Modular counting quantifiers}\label{ss:rmc}

We briefly recall the definition of a modular counting quantifier. Given $\phi(\bar x)=\C^{r,m}_x\psi(\bar x\cup\{x\})$, this expression is said to hold in $\mathfrak M$ for a specific assignment to the variable of $\bar x$, if the size of the set $\{x:\mathfrak M\models\psi(\bar x\cup\{x\})\}$ is congruent to $r$ modulo $m$. As with the previous subsection, we assume that the quantified variable is not used outside the quantification scope, and that no variable is quantified more than once. We again denote for notational convenience the quantified variable by $x=x_{\ell_{\bar x}+1}$, and note that $\mathfrak X\subseteq [\ell_{\bar x}]$ cannot include a reference to $x$.

When working with $(\mathfrak N,\mathfrak U,\mathfrak R)$ that corresponds to $\mathfrak M$, to obtain the original modular count, we have to count the set (satisfying the induction hypothesis) $\{x:\mathfrak N\models\psi'_{\mathfrak X,\mathfrak U,\mathfrak R}(\bar x\cup\{x\}\setminus x_{\mathfrak X})\}$, as well as check whether $\mathfrak N\models\psi'_{\mathfrak X\cup\{\ell_{\bar x}+1\},\mathfrak U,\mathfrak R}(\bar x\setminus x_{\mathfrak X})$ (which if true adds $1$ to the count). This gives $(\C^{r-1,m}_x\!\psi'_{\mathfrak X,\mathfrak U,\mathfrak R}(\bar x\cup\{x\}\!\setminus\! x_{\mathfrak X})\wedge\psi'_{\mathfrak X\cup\{\ell_{\bar x}+1\},\mathfrak U,\mathfrak R}(\bar x\setminus\! x_{\mathfrak X}))\vee(\C^{r,m}_x\!\psi'_{\mathfrak X,\mathfrak U,\mathfrak R}(\bar x\cup\{x\}\!\setminus\! x_{\mathfrak X})\wedge\neg\psi'_{\mathfrak X\cup\{\ell_{\bar x}+1\},\mathfrak U,\mathfrak R}(\bar x\setminus\! x_{\mathfrak X}))$ as the combined expression for $\phi'_{\mathfrak X,\mathfrak U,\mathfrak R}(\bar x\setminus x_{\mathfrak X})$.

\subsection{Monadic second order quantifiers}\label{ss:rmso}

Here we deal with quantifiers over unary relations. The cases we cover are the existential quantification $\phi(\bar x)=\exists_U\psi(\bar x)$ and the universal quantification $\phi(\bar x)=\forall_U\psi(\bar x)$, where $U$ is a new unary relation that does not appear in the language $(\bar a,\bar U,\bar R)$ of $\phi(\bar x)$, while being part of the language of $\psi(\bar x)$. As before, we assume that the quantified relation symbol $U$ appears only in the scope of this quantification, and is not quantified anywhere else, and again denote for convenience $U=U_{\ell_{\bar U}+1}$. In particular, when analyzing expressions of the type $\psi'_{\mathfrak X,\mathfrak U',\mathfrak R}(\bar x\setminus x_{\mathfrak X})$, we may allow $\mathfrak U'$ to contain $[\ell_{\bar U}+1]$ (the same is not allowed for the expression $\phi'_{\mathfrak X,\mathfrak U,\mathfrak R}(\bar x\setminus x_{\mathfrak X})$, whose language does not contain $U$).

Consider now the family of possible models $\mathfrak M'$ that extend $\mathfrak M$ with an interpretation of the relation $U$. Now consider $(\mathfrak N',\mathfrak U',\mathfrak R')$ which correspond to $\mathfrak M'$, in relation to $(\mathfrak N,\mathfrak U,\mathfrak R)$ which correspond to $\mathfrak M$. Referring to Definition \ref{df:corr}, every relation already appearing in $\bar U$ will have the same interpretation in $\mathfrak N$ and $\mathfrak N'$. Also, $\mathfrak R'=\mathfrak R$, since the binary relations are the same in the languages of both models. Additionally, from the definition, the interpretation of $U=U_{\ell_{\bar U}+1}$ in $\mathfrak N'$ is the restriction of its interpretation in $\mathfrak M'$ to $[n+\ell_{\bar a}-1]$. As for the final ingredient $\mathfrak U'$, for every $i\in [\ell_{\bar U}]$, the condition on whether it is in $\mathfrak U$ or in $\mathfrak U'$ is the same. However, $\mathfrak U'$ may also include $\ell_{\bar U}+1$ according to whether $\mathfrak M'\models U(a)$. So considering all possible models $\mathfrak M'$, we have two possibilities. Either $\mathfrak U'=\mathfrak U$, or $\mathfrak U'=\mathfrak U\cup\{\ell_{\bar U}+1\}$.

We are now ready to construct our expression that corresponds to all models extending $\mathfrak M$. For the existential case we have $\phi'_{\mathfrak X,\mathfrak U,\mathfrak R}(\bar x\setminus x_{\mathfrak X})=\exists_U\psi'_{\mathfrak X,\mathfrak U,\mathfrak R}(\bar x\setminus x_{\mathfrak X})\vee\exists_U\psi'_{\mathfrak X,\mathfrak U\cup\{\ell_{\bar U}+1\},\mathfrak R}(\bar x\setminus x_{\mathfrak X})$, and for the universal one we have $\phi'_{\mathfrak X,\mathfrak U,\mathfrak R}(\bar x\setminus x_{\mathfrak X})=\forall_U\psi'_{\mathfrak X,\mathfrak U,\mathfrak R}(\bar x\setminus x_{\mathfrak X})\wedge\forall_U\psi'_{\mathfrak X,\mathfrak U\cup\{\ell_{\bar U}+1\},\mathfrak R}(\bar x\setminus x_{\mathfrak X})$.

\section{Nullary relations and a many-one version of the reduction}\label{s:m}

Nullary relations are relations of arity zero. Formally, for a nullary relation $Z$, the corresponding atomic formula is $Z()$, and a model $\mathfrak M$ over a language that includes $Z$ interprets this formula as either true or false, that is, either $\mathfrak M\models Z()$ or $\mathfrak M\models\neg Z()$.

Note that nullary relations can be simulated using higher arity relations. To replace a nullary relation $Z$ in the language with a unary relation $U$ (while preserving the model count), the logical expression under discussion should be conjuncted with ``$\forall_x\forall_y(U(x)\leftrightarrow U(y))$'', and then every instance of ``$Z()$'' in the formula should be replaced with ``$\exists_xU(x)$''.

For convenience, in the following we use explicit nullary relations in our formalism. We prove in this section a ``many-one'' reduction of counting with hard-coded constants to counting without them, that is, a reduction of the counting function to another counting function based on a single class, rather than a reduction to the sum of several such functions.

\begin{theorem}[Many-one reduction to the case without constants]\label{th:mr}
For any class $\mathcal C$ defined by an FOL (resp.\ MSOL, CMSOL) sentence involving a set of constant symbols $\bar a$, nullary symbols $\bar Z$, unary symbols $\bar U$ and binary symbols $\bar R$, there exists a class $\mathcal C'$ definable by an FOL (resp.\ MSOL, CMSOL) sentence involving $\bar Z'$, $\bar U'$ (which contain $\bar Z$ and $\bar U$ respectively), $\bar R$ and no constants, satisfying $f_{\mathcal C}(n)=f_{\mathcal C'}(n)$ for all $n\in\mathbb N$.
\end{theorem}

Also here, the theorem follows from a single constant removal lemma, which is used for an inductive argument over $\ell_{\bar a}$.

\begin{lemma}[Removing a single constant in a many-one manner]\label{lm:moscr}
For any class $\mathcal C$ defined by an FOL (resp.\ MSOL, CMSOL) sentence involving a set of constant symbols $\bar a$ with $\ell_{\bar a}>0$, nullary symbols $\bar Z$, unary symbols $\bar U$ and binary symbols $\bar R$, there exists a class $\mathcal C'$, definable by an FOL (resp.\ MSOL, CMSOL) sentence over the language $(\bar a',\bar Z',\bar U',\bar R')$, where $\bar a'=\bar a\setminus\{a_{\ell_{\bar a}}\}$, $\bar Z'=\bar Z\cup\bar S\cup\bar D$ where $\ell_S=\ell_U$ and $\ell_D=\ell_R$, $\bar U'=\bar U\cup\bar I\cup\bar O$ where $\ell_{\bar I}=\ell_{\bar O}=\ell_{\bar R}$, and $\bar R'=\bar R$, satisfying $f_{\mathcal C}(n)=f_{\mathcal C'}(n)$ for all $n\in\mathbb N$.
\end{lemma}

The main new idea in the proof of this version is to use new nullary relations to hold the information as to whether $R(a,a)$ holds for a binary relation $R$, or whether $U(a)$ holds for a unary relation $U$, while in the original version we constructed different expressions for each of these options. So given $\phi(\bar x)$, our inductive construction will produce the expression $\phi'_{\mathfrak X}(\bar x\setminus x_{\mathfrak X})$ for every possible $\mathfrak X\subseteq [\ell_{\bar x}]$, without referring to the sets $\mathfrak U$ and $\mathfrak R$ that appeared in the proof of Lemma \ref{lm:scr} and held the information about which relations contain $a$ or $(a,a)$.

The definition of correspondence is adapted to use the new nullary relations in the language, instead of referring to prescribed sets, in the following way.

\begin{definition}\label{df:mocorr}
Given a model $\mathfrak M$ over the language $(\bar a,\bar Z,\bar U,\bar R)$ with universe $[n+\ell_{\bar a}]$, and a model $\mathfrak N$ over the language $(\bar a',\bar Z',\bar U',\bar R)$ with universe $[n+\ell_{\bar a}-1]$, where in both models every constant $a_i$ is interpreted to be $n+i$, we say that $\mathfrak N$ {\em corresponds} to $\mathfrak M$ if the following hold.
\begin{itemize}
\item For every $Z\in \bar Z$, we have $\mathfrak N\models Z()$ if and only if $\mathfrak M\models Z()$.
\item For every $U\in \bar U$ and $x\in [n+\ell_{\bar a}-1]$, we have $\mathfrak N\models U(x)$ if and only if $\mathfrak M\models U(x)$.
\item For every $i\in [\ell_{\bar U}]$, we have $\mathfrak N\models S_i()$ if and only if $\mathfrak M\models U_i(a)$.
\item For every $R\in \bar R$ and $x,y\in [n+\ell_{\bar a}-1]$, we have $\mathfrak N\models R(x,y)$ if and only if $\mathfrak M\models R(x,y)$.
\item For every $i\in [\ell_{\bar R}]$ and $x\in [n+\ell_{\bar a}-1]$, we have $\mathfrak N\models I_i(x)$ if and only if $\mathfrak M\models R_i(x,a)$.
\item For every $i\in [\ell_{\bar R}]$ and $x\in [n+\ell_{\bar a}-1]$, we have $\mathfrak N\models O_i(x)$ if and only if $\mathfrak M\models R_i(a,x)$.
\item For every $i\in [\ell_{\bar R}]$, we have $\mathfrak N\models D_i()$ if and only if $\mathfrak M\models R_i(a,a)$.
\end{itemize}
\end{definition}

As expected we have the immediate counterpart to Observation \ref{ob:bij}.

\begin{observation}\label{ob:mobij}
Definition \ref{df:mocorr} provides a bijection between the set of possible models $\mathfrak M$ over the language $(\bar a,\bar Z,\bar U,\bar R)$ with universe $[n+\ell_{\bar a}]$, and the set of possible models $\mathfrak N$ over the language $(\bar a',\bar Z',\bar U',\bar R)$ with universe $[n+\ell_{\bar a}-1]$ (where all constants are interpreted as in Definition \ref{df:mocorr}).
\end{observation}

The rest of this section is concerned with the recursive definition of $\phi'_{\mathfrak X}(\bar x\setminus x_{\mathfrak X})$ given $\phi(\bar x)$, satisfying the following correspondence invariant: If $\mathfrak M$ and $\mathfrak N$ correspond according to Definition \ref{df:mocorr}, then for any fixing of the variables $\bar x\setminus x_{\mathfrak X}$, we have $\mathfrak M\models\phi(\bar x_{\mathfrak X\to a})$ if and only if $\mathfrak N\models\phi'_{\mathfrak X}(\bar x\setminus x_{\mathfrak X})$ (we refer to the Section \ref{s:r} for the definitions of $\bar x\setminus x_{\mathfrak X}$ and $\bar x_{\mathfrak X\to a}$). The main differences between this construction and the one of Section \ref{s:r} are in the handling of atomic formulas and of monadic second order quantifiers.

\subsection{The base constructions}\label{ss:morbase}

For the Boolean ``true'' and ``false'' statements, just as with Subsection \ref{ss:rbase}, if $\phi(\bar x)$ is $\top$ (respectively $\bot$), then we also set $\phi'_{\mathfrak X}(\bar x\setminus x_{\mathfrak X})$ to $\top$ (respectively $\bot$).

For the statement $\phi(\bar x)=Z_i()$ where $i\in [\ell_{\bar Z}]$, relating to a nullary relation that was already present in the language of $\mathfrak M$, we simply set $\phi'_{\mathfrak X}(\bar x\setminus x_{\mathfrak X})=Z_i()$ as well.

For the expression $\phi(\bar x)=U_i(x_j)$ where $i\in [\ell_{\bar U}]$ and $j\in [\ell_{\bar x}]$, to produce $\phi'_{\mathfrak X}(\bar x\setminus x_{\mathfrak X})$, we again partition to cases according to whether $j\in\mathfrak X$. In the case where $j\notin\mathfrak X$, we again simply set $\phi'_{\mathfrak X}(\bar x\setminus x_{\mathfrak X})=U_i(x_j)$ to satisfy the correspondence invariant (recall that the ``if and only if'' condition in this case should hold when the value of $x_i$ is in $[n+\ell_{\bar a}-1]$).

Similarly, for $i\in [\ell_{\bar U}]$ and $j\in [\ell_{\bar a}-1]$, for $\phi(\bar x)=U_i(a_j)$ we set $\phi'_{\mathfrak X}(\bar x\setminus x_{\mathfrak X})=U_i(a_j)$, noting that in our setting the value of $a_j$ is guaranteed to equal $n+j\in [n+\ell_{\bar a}-1]$.

Returning to the expression $\phi(\bar x)=U_i(x_j)$, for the case where $x_j\in\mathfrak X$, recall that here $\phi'_{\mathfrak X}(\bar x\setminus x_{\mathfrak X})$ should not depend on $x_j$, but should rather express (for the correspondence invariant to hold) whether $\mathfrak M\models U_i(a)$. For the corresponding model $\mathfrak N$ this information is kept by the nullary relation $S_i$. Therefor, for the invariant to hold, we set $\phi'_{\mathfrak X}(\bar x\setminus x_{\mathfrak X})$ to $S_i()$.

The remaining case for a unary relation is the expression $\phi(\bar x)=U_i(a)$. Also here we set $\phi'_{\mathfrak X}(\bar x\setminus x_{\mathfrak X})$ to $S_i()$.

Moving on to the atomic expressions involving a binary relation $R_i$ where $i\in [\ell_{\bar R}]$, the partition to cases is analogous to that of Subsection \ref{ss:rbase}. The first case is where $\phi(\bar x)=R_i(x_j,x_k)$ where both $j\notin\mathfrak X$ and $k\notin\mathfrak X$, for which as expected we set $\phi'_{\mathfrak X}(\bar x\setminus x_{\mathfrak X})=R_i(x_j,x_k)$, using very much the same argument.

The next four cases are those that resemble the last two cases that we discussed about a unary relation. Namely, these are the cases where $\phi(\bar x)=R_i(x_j,x_k)$ with $j,k\in\mathfrak X$, $\phi(\bar x)=R_i(x_j,a)$ or $\phi(\bar x)=R_i(a,x_j)$ with $j\in\mathfrak X$, and $\phi(\bar x)=R_i(a,a)$. In all these cases the resulting expression should reflect on whether $\mathfrak M\models R_i(a,a)$, which for the corresponding $\mathfrak N$ is handled by the nullary relation $D_i$. We hence set $\phi'_{\mathfrak X}(\bar x\setminus x_{\mathfrak X})=D_i()$.

The remaining cases are handled identically to Subsection \ref{ss:rbase}. For the case $\phi(\bar x)=R_i(x_j,x_k)$ where $j\notin\mathfrak X$ and $k\in\mathfrak X$, and the case $\phi(\bar x)=R_i(x_j,a)$ where $j\notin\mathfrak X$, for the correspondence invariant to follow we set $\phi'_{\mathfrak X}(\bar x\setminus x_{\mathfrak X})=I_i(x_j)$. For the similar cases $\phi(\bar x)=R_i(a_j,x_k)$ and $\phi(\bar x)=R_i(a_j,a)$, where $j\in[\ell_{\bar a}-1]$ and $k\in\mathfrak X$, we set $\phi'_{\mathfrak X}(\bar x\setminus x_{\mathfrak X})=I_i(a_j)$. For the case $\phi(\bar x)=R_i(x_k,x_j)$ where $j\notin\mathfrak X$ and $k\in\mathfrak X$, and the case $\phi(\bar x)=R_i(a,x_j)$ where $j\notin\mathfrak X$, we set $\phi'_{\mathfrak X}(\bar x\setminus x_{\mathfrak X})=O_i(x_j)$. Finally, for the cases $\phi(\bar x)=R_i(x_k,a_j)$ and $\phi(\bar x)=R_i(a,a_j)$ where $j\in[\ell_{\bar a}-1]$ and $k\in\mathfrak X$, we set $\phi'_{\mathfrak X}(\bar x\setminus x_{\mathfrak X})=O_i(a_j)$.

The conversion of atomic expressions using the equality relation is completely identical to that of Section \ref{ss:rbase}.

\subsection{Boolean connectives}\label{ss:morbool}

The handling of Boolean connectives is completely identical to that of Subsection \ref{ss:rbool}, including the same straightforward arguments. For example, for the expression $\phi(\bar x)=\psi_1(\bar x)\wedge\psi_2(\bar x)$ we set $\phi'_{\mathfrak X}(\bar x\setminus x_{\mathfrak X})=\psi'_{1,\mathfrak X}(\bar x\setminus x_{\mathfrak X})\wedge\psi'_{2,\mathfrak X}(\bar x\setminus x_{\mathfrak X})$.

\subsection{Quantifiers over variables}\label{ss:morfocm}

The handling of first order quantifiers, as well as modular counting quantifiers, is pretty much the same as that of the respective Subsection \ref{ss:rfo} and Subsection \ref{ss:rmc}. In particular, we again require that every quantified variable is quantified only once in the expression, and is not used at all outside the scope of the quantification, meaning in particular that the set $\mathfrak X$ cannot contain a reference to the quantified variable.

For notational convenience, when $\phi(\bar x)$ is our formula, we again denote by $x=x_{\ell_{\bar x}+1}$ the quantified variable, so for example the existential quantification case is $\phi(\bar x)=\exists_x\psi(\bar x\cup\{x\})$. For each case we will construct a corresponding $\phi'_{\mathfrak X}(\bar x\setminus x_{\mathfrak X})$ where $\mathfrak X\subseteq [\ell_{\bar x}]$.

In the existential case, we want $\mathfrak N\models\phi'_{\mathfrak X}(\bar x\setminus x_{\mathfrak X})$ to occur whenever there is at least one value of $x$ for which $\mathfrak M\models\psi(\bar x\cup\{x\})$. If this occurs for a value of $x$ within $[n+\ell_{\bar a}-1]$, then by the induction hypothesis this is covered by the expression $\exists_x\psi'_{\mathfrak X}(\bar x\cup\{x\}\setminus x_{\mathfrak X})$, and if the above occurs in $\mathfrak M$ for $x=a$, then this is covered by the expression $\psi'_{\mathfrak X\cup\{\ell_{\bar x}+1\}}(\bar x\setminus x_{\mathfrak X})$. The combined expression that satisfies the correspondence invariant is hence $\phi'_{\mathfrak X}(\bar x\setminus x_{\mathfrak X})=\exists_x\psi'_{\mathfrak X}(\bar x\cup\{x\}\setminus x_{\mathfrak X})\vee\psi'_{\mathfrak X\cup\{\ell_{\bar x}+1\}}(\bar x\setminus x_{\mathfrak X})$.

The universal case follows an analogous argument, only here $\mathfrak M\models\psi(\bar x\cup\{x\})$ needs to hold for all values of $x$, those in $[n+\ell_{\bar a}-1]$ as well as the value $a$. The combined expression is hence $\phi'_{\mathfrak X}(\bar x\setminus x_{\mathfrak X})=\forall_x\psi'_{\mathfrak X}(\bar x\cup\{x\}\setminus x_{\mathfrak X})\wedge\psi'_{\mathfrak X\cup\{\ell_{\bar x}+1\}}(\bar x\setminus x_{\mathfrak X})$.

We now handle the modular counting quantifier case $\phi(\bar x)=\C^{r,m}_x\psi(\bar x\cup\{x\})$. To obtain the original modular count for $\mathfrak M$ when working with the corresponding $\mathfrak N$, we have to count the set (satisfying the induction hypothesis) $\{x:\mathfrak N\models\psi'_{\mathfrak X}(\bar x\cup\{x\}\setminus x_{\mathfrak X})\}$, as well as check whether $\mathfrak N\models\psi'_{\mathfrak X\cup\{\ell_{\bar x}+1\}}(\bar x\setminus x_{\mathfrak X})$. For the complete combined expression for $\phi'_{\mathfrak X}(\bar x\setminus x_{\mathfrak X})$ we obtain $(\C^{r-1,m}_x\psi'_{\mathfrak X}(\bar x\cup\{x\}\setminus x_{\mathfrak X})\wedge\psi'_{\mathfrak X\cup\{\ell_{\bar x}+1\}}(\bar x\setminus x_{\mathfrak X}))\vee(\C^{r,m}_x\psi'_{\mathfrak X}(\bar x\cup\{x\}\setminus x_{\mathfrak X})\wedge\neg\psi'_{\mathfrak X\cup\{\ell_{\bar x}+1\}}(\bar x\setminus x_{\mathfrak X}))$.

\subsection{Monadic second order quantifiers}\label{ss:mormso}

We now deal with the  cases of existential quantification $\phi(\bar x)=\exists_U\psi(\bar x)$ and universal quantification $\phi(\bar x)=\forall_U\psi(\bar x)$, where $U$ is a new unary relation that does not appear in the language $(\bar a,\bar Z,\bar U,\bar R)$ of $\phi(\bar x)$, while being part of the language of $\psi(\bar x)$. As before, we assume that the quantified relation symbol $U$ appears only in the scope of this quantification, and is not quantified anywhere else, and denote for convenience $U=U_{\ell_{\bar U}+1}$.

When preparing to analyze expressions of the type $\psi'_{\mathfrak X}(\bar x\setminus x_{\mathfrak X})$, as per the induction hypothesis for the construction of $\phi'_{\mathfrak X}(\bar x\setminus x_{\mathfrak X})$, we note that just as $U_1,\ldots,U_{\ell_{\bar U}}$ have their counterpart nullary relations $S_1,\ldots,S_{\ell_{\bar U}}$, we also need a new counterpart $S=S_{\ell_{\bar U}+1}$ to $U=U_{\ell_{\bar U}+1}$. When analyzing a model $\mathfrak M'$ for $\psi(\bar x)$, which (unlike $\mathfrak M$) interprets $U$ as well, we note that the corresponding $\mathfrak N'$ must interpret both $U$ and $S$, where in particular $\mathfrak N'\models S()$ if and only if $\mathfrak M'\models U(a)$. When constructing $\phi'_{\mathfrak X}(\bar x\setminus x_{\mathfrak X})$, we will need to quantify both $U$ and $S$. Quantifying over $S$ as well as $U$ makes sure that our quantification encompasses, for the tests against the original $\psi(\bar x)$, both $U$ for which $U(a)$ holds and $U$ for which $U(a)$ does not hold. Since nullary relations can be simulated by unary relations, utilizing the notion of quantification over nullary relations does not move us outside the realm of monadic second order logic.

Having discussed the role of the new relation, constructing the expressions that relate to all models corresponding to those extending $\mathfrak M$ is now simple. For the existential case we set $\phi'_{\mathfrak X}(\bar x\setminus x_{\mathfrak X})=\exists_U\exists_S\psi'_{\mathfrak X}(\bar x\setminus x_{\mathfrak X})$, and for the universal case we set $\phi'_{\mathfrak X}(\bar x\setminus x_{\mathfrak X})=\forall_U\forall_S\psi'_{\mathfrak X}(\bar x\setminus x_{\mathfrak X})$.

\section{Higher arity relations}\label{s:h}

We prove in this section an extension of Theorem \ref{th:mr} to higher arity relations.

\begin{theorem}[Many-one reduction allowing higher arity]\label{th:hr}
For any class $\mathcal C$ defined by an FOL (resp.\ MSOL, CMSOL, GSOL, SOL) sentence involving a set of constant symbols $\bar a$, and relation symbols $\bar R$ of arbitrary arities, there exists a class $\mathcal C'$ definable by an FOL (resp.\ MSOL, CMSOL, GSOL, SOL) sentence involving $\bar R'$, which contains $\bar R$, has the same maximum arity as $\bar R$, and has no new relations of maximum arity, satisfying $f_{\mathcal C}(n)=f_{\mathcal C'}(n)$ for all $n\in\mathbb N$.
\end{theorem}

Also here, this follows by induction of the following lemma, extending Lemma \ref{lm:moscr}

\begin{lemma}[Removing a single constant with higher arities]\label{lm:hscr}
For any class $\mathcal C$ defined by an FOL (resp.\ MSOL, CMSOL, GSOL, SOL) sentence involving a set of constant symbols $\bar a$, and relation symbols $\bar R$ of arbitrary arities, there exists a class $\mathcal C'$, definable by an FOL (resp.\ MSOL, CMSOL, GSOL, SOL) sentence over the language $(\bar a',\bar R')$, where $\bar a'=\bar a\setminus\{a_{\ell_{\bar a}}\}$, $\bar R'$ contains $\bar R$, has the same maximum arity as $\bar R$, and has no new relations of maximum arity, satisfying $f_{\mathcal C}(n)=f_{\mathcal C'}(n)$ for all $n\in\mathbb N$.
\end{lemma}

The construction is very similar to that of Section \ref{s:m}, and we only highlight here the differences, which are mainly in the definition of the corresponding models and the transformation of atomic relational formulas.

We first set up the language: We assume that $\bar R=\{R_1,\ldots,R_{\ell_{\bar R}}\}$ are relation symbols whose arities are $\bar\rho=\{\rho_1,\ldots,\rho_{\ell_{\bar R}}\}$ respectively. To construct $\bar R'$, first every relation $R_i\in\bar R$ is replaced with $2^{\rho_i}$ relation symbols $\bar R_i=\langle R_{i,A}:A\subseteq [\rho_i]\rangle$, where each relation $R_{i,A}$ is of arity $\rho_i-|A|$. Note that in particular $R_{i,\emptyset}$ is identified with the original $R_i$, and that $R_{i,[\rho_i]}$ is a nullary relation. Finally, $\bar R'$ is the union of these sets, $\bar R'=\bigcup_{i=1}^{\ell_{\bar R}}\bar R_i$.

For some intuition of why this language is a straightforward extension of the one defined in conjunction with Definition \ref{df:mocorr} in Section \ref{s:m}, consider the case of a binary relation $R_i$. In this case, $R_{i,\{1\}}$ is the same as $O_i$ in Section \ref{s:m}, $R_{i,\{2\}}$ is the same as $I_i$ there, and $R_{i,\{1,2\}}$ is the same as $D_i$ there. The following is a generalization of Definition \ref{df:mocorr}.

\begin{definition}\label{df:hcorr}
Given a model $\mathfrak M$ over the language $(\bar a,\bar R)$ with universe $[n+\ell_{\bar a}]$, and a model $\mathfrak N$ over the language $(\bar a',\bar R')$ with universe $[n+\ell_{\bar a}-1]$, where in both models every constant $a_i$ is interpreted to be $n+i$, we say that $\mathfrak N$ {\em corresponds} to $\mathfrak M$ if the following holds.
\begin{itemize}
\item For every $R_i\in\bar R$ and every $\bar x=x_1,\ldots,x_{\rho_i}\in[n+\ell_{\bar a}]$, denoting by $A$ the set of indexes of the variables for which $x=a=a_{\ell_{\bar a}}$, that is $A=\{i:x_i=a\}$, we have $\mathfrak M\models R(x_1,\ldots,x_{\rho_i})$ if and only if $\mathfrak N\models R_{i,A}(\bar x\setminus x_A)$, where as before we let $\bar x\setminus x_A$ denote the subsequence of variables whose indexes are not in $A$.
\end{itemize}
\end{definition}

The exact analog to Observation \ref{ob:mobij} also holds, and we define the same correspondence invariant with respect to logic expressions. The rest of this section is concerned with the recursive definition of $\phi'_{\mathfrak X}(\bar x\setminus x_{\mathfrak X})$ given $\phi(\bar x)$ and $\mathfrak X\subseteq [\ell_{\bar x}]$, satisfying the correspondence invariant.

\subsection{FOL expressions and counting quantifiers}

The only difference between this section and Subsections \ref{ss:morbase}, \ref{ss:morbool} and \ref{ss:morfocm} is in the construction of $\phi'_{\mathfrak X}(\bar x\setminus x_{\mathfrak X})$ where $\phi(\bar x)=R_i(y_1,\ldots,y_{\rho_i})$, and each $y_j$ is either some variable $x_{i_j}$ or some constant $a_{i_j}$.

For this construction, we let the set $A$ denote the indexes of all $j$ for which $y_j$ is either some $x_i$ for $i\in\mathfrak X$, or the constant $a$ (but not any constant $a_{i_j}$ for $i_j<\ell_{\bar a}$). We denote by $\bar y\setminus y_A$ the subsequence of $\bar y=y_1,\ldots,y_{\rho_i}$ after excluding all $y_j$ with $j\in A$, and define $\phi'_{\mathfrak X}(\bar x\setminus x_{\mathfrak X})=R_{i,A}(\bar y\setminus y_A)$.

All other base constructions, as well as the recursive constructions for Boolean connective and quantification over variables, are identical to those of Section \ref{s:m}.

\subsection{Second order quantifiers}

We first look at the case of existential quantification $\phi(\bar x)=\exists_R\psi(\bar x)$, where $R$ is a new relation that does not appear in the language $(\bar a,\bar R)$ of $\phi(\bar x)$, while being part of the language of $\psi(\bar x)$. Again we assume that $R$ appears only in the scope of this quantification, and is not quantified anywhere else, and denote for convenience $R=R_{\ell_{\bar R}+1}$. We also denote by $\rho=\rho_{\ell_{\bar R}+1}$ the arity of $R$.

To construct $\phi'_{\mathfrak X}(\bar x\setminus x_{\mathfrak X})$, just as we expanded each $R_i$ to $2^{\rho_i}$ relations for $1\leq i\leq \ell_{\bar R}$, we expand the quantified relation $R$ to a sequence of $2^{\rho}$ relations $\langle R_A:A\subseteq [\rho]\rangle$, where each $R_A$ is of arity $\rho-|A|$ ($R_A$ is not part of the relations of the language of $\phi'_{\mathfrak X}(\bar x\setminus x_{\mathfrak X})$, but identified with $R_{\ell_{\bar R}+1,A}$ it is part of the language of $\psi'_{\mathfrak X}(\bar x\setminus x_{\mathfrak X})$ that is constructed by induction from $\psi(\bar x)$).

The quantification will be over all new relations, that is $\phi'_{\mathfrak X}(\bar x\setminus x_{\mathfrak X})=\exists_{\langle R_A:A\subseteq [\rho]\rangle}\psi'_{\mathfrak X}(\bar x\setminus x_{\mathfrak X})$. Note that in particular for an MSOL quantification, that is when $R$ is of arity $1$, we will have a quantification over a unary relation $R_{\emptyset}$ and over a nullary relation $R_{\{1\}}$, just as with the construction in Subsection \ref{ss:mormso}.

For the case of universal quantification $\phi(\bar x)=\forall_R\psi(\bar x)$ we again define $\langle R_A:A\subseteq [\rho]\rangle$, and analogously set $\phi'_{\mathfrak X}(\bar x\setminus x_{\mathfrak X})=\forall_{\langle R_A:A\subseteq [\rho]\rangle}\psi'_{\mathfrak X}(\bar x\setminus x_{\mathfrak X})$.

Finally, we briefly explain how to deal with guarded second order (GSOL) quantifiers. These are quantifiers over a new relation $R$ whose arity is identical to that of an existing relation $R_i$, where we look only at the possibilities for $R$ that make it a subset of $R_i$. The existential case is written $\phi(\bar x)=\exists_{R\subseteq R_i}\psi(\bar x)$, and the universal case is written $\phi(\bar x)=\forall_{R\subseteq R_i}\psi(\bar x)$.

We construct the relations $\langle R_A:A\subseteq [\rho]\rangle$, where the arity $\rho_i-|A|$ of $R_A$ is identical to that of $R_{i,A}$ in the language of $\phi'_{\mathfrak X}(\bar x\setminus x_{\mathfrak X})$. We simply quantify every $R_A$ as a subset of its respective $R_{i,A}$, so for the existential case we have $\phi'_{\mathfrak X}(\bar x\setminus x_{\mathfrak X})=\exists_{\langle R_A\subseteq R_{i,A}:A\subseteq [\rho]\rangle}\psi'_{\mathfrak X}(\bar x\setminus x_{\mathfrak X})$, and for the universal case we have $\phi'_{\mathfrak X}(\bar x\setminus x_{\mathfrak X})=\forall_{\langle R_A\subseteq R_{i,A}:A\subseteq [\rho]\rangle}\psi'_{\mathfrak X}(\bar x\setminus x_{\mathfrak X})$.

\part{\Large Ternary relations}
\label{part-3}
\section{An $\FOL$-definable class $\cC$ where $f_{\cC}(n)$ is not MC-finite}
\label{se:ternary}

In this part we negatively settle the question of whether the Specker-Blatter theorem holds for classes whose 
language contains only ternary and lower-arity relations. 

\subsection{Using one hard-wired constant}
We first construct a class whose language includes a single ternary relation and a single hard-coded constant.
Our counterexample builds on ideas used in \cite{noquad}.

\begin{theorem}[Ternary relation counterexample with a constant]
\label{th:tce}
There exists an FOL sentence $\phi_{\mathcal M}$ over the language $(a,R)$, where $a$ is a single (hard-coded) constant and $R$ is a single relation of arity $3$, so that the corresponding class $\mathcal C$ satisfies $f_{\mathcal C}(n-1)=0$ for any $n$ that is not a power of $2$, and $f_{\mathcal C}(n-1)\equiv 1\pmod{2}$ for $n=2^m$ for every $m\in\mathbb{N}$. In particular, $f_{\mathcal C}$ is not ultimately periodic modulo $2$.
\end{theorem}

The statement uses $f_{\mathcal C}(n-1)$ instead of $f_{\mathcal C}(n)$, but recalling the definition of $f_{\mathcal C}$, this refers to the universe $[n-1]\cup\{a\}$ whose size is $n$. We explain later how to modify this class to produce a counterexample modulo other prime numbers $p$ instead of $2$.

By Theorem \ref{th:hr}, we have the following immediate corollary that does away with the constant, at the price of adding some additional smaller arity relations.

\begin{corollary}[Ternary counterexample without constants]\label{cor:ptce}
There exists an FOL sentence $\phi'_{\mathcal M}$ over the language $(\bar R)$, where $\bar R$ includes one relation of arity $3$ and other relations of lower arities, so that the corresponding class $\mathcal C$ satisfies $f_{\mathcal C}(n)=0$ for every $n$ for which $n+1$ is not a power of $2$, and $f_{\mathcal C}(n)\equiv 1\pmod{2}$ for $n=2^m-1$ for every $m\in\mathbb{N}$. In particular, $f_{\mathcal C}$ is not ultimately periodic modulo $2$.
\end{corollary}

In Section \ref{s:s} we will show how to further reduce the language so that it includes only one ternary relation and no lower arity relations.

\subsection{The first construction}
The starting point of the construction is a structure that is defined over a non-constant length sequence (and hence not yet expressible in FOL) of unordered graphs. This definition follows the streamlining by Specker \cite{spadd} of the original construction from \cite{noquad}.

\begin{definition}[Iterated matching sequence]
Given a set $V$ of vertices, An {\em iterated matching sequence} is a sequence of graphs over $V$, identified by their edge sets $\bar E=E_1,\ldots,E_{\ell_{\bar E}}$, satisfying the following for every $1\leq i\leq \ell_{\bar E}$.
\begin{itemize}
\item The connected components of $E_i$ are (vertex-disjoint) complete bipartite graphs.
\item The two vertex classes of every complete bipartite graph in $E_i$ as above are two connected components of $\bigcup_{j=1}^{i-1}E_j$ (for $i=1$ this means that $E_1$ is a matching).
\item Every connected component of $\bigcup_{j=1}^{i-1}E_j$ is a vertex class of some bipartite graph of $E_i$ (so in particular $E_1$ is a perfect matching).
\end{itemize}
An iterated matching sequence $\bar E$ is {\em full} if every vertex pair $u,v\in V$ (where $u\neq v$) appears in some $E_i$.
\end{definition}

The following properties of iterated matching sequences are easily provable by induction.

\begin{observation}\label{ob:im}
For an iterated matching $\bar E$, every $E_i$ corresponds to a perfect matching over the set of connected components of $\bigcup_{j=1}^{i-1}E_j$. Additionally, every connected component of $\bigcup_{j=1}^iE_j$ is a clique with exactly $2^i$ vertices.
\end{observation}

The above implies that there can be a full iterated matching sequence over $[n]$ if and only if $n$ is a power of $2$, in which case $\ell_{\bar E}=\log_2(n)$. Denoting the number of possible full iterated matching sequences over $[n]$ by $f_{\mathcal M}(n)$, note the following lemma.

\begin{lemma}[see \cite{spadd}]
For every $n$ which is not a power of $2$ we have $f_{\mathcal M}(n)=0$, while $f_{\mathcal M}(n)\equiv 1\pmod{2}$ for $n=2^m$ for every $m\in\mathbb N$.
\end{lemma}

The rest of this section concerns the construction of a sentence $\phi_{\mathcal M}$ over a language with one constant and one ternary relation, so that the corresponding class $\mathcal C$ satisfies $f_{\mathcal C}(n-1)=f_{\mathcal M}(n)$. In the original construction utilizing a quaternary relation $Q$, essentially we had $(u,v,x,y)\in Q$ if $(u,v)\in E_i$ and $(x,y)\in E_{i-1}$ for some $1<i\leq\ell_{\bar E}$, or $(u,v)\in E_1$ and $x=y$. For the construction here, we only have a ternary relation $R$, and we encode the placement of $(u,v)$ within $\bar E$ by the set $\{w:(u,v,w)\in R\}$. We will have to utilize the hard-coded constant $a$ to make sure that there is exactly one way to encode every full iterated matching sequence.

\subsection{Setting up and referring to an order over the vertex pairs}\label{ss:to}

We simulate the structure of a full iterated matching sequence over $[n]$ (where $n\in [n]$ is identified with the constant $a$) by assigning ``ranks'' to pairs of members of $[n]$, which we consider as vertices, where each pair $(x,y)$ is assigned the set $r_{x,y}=\{z:(x,y,z)\in R\}$. First we need to make sure that ``graphness'' is satisfied, which means that $r_{x,y}$ is symmetric and is empty for loops.

$$\phi_{\mathrm{graph}}=\forall_{x,y,z}(R(x,y,z)\to(x\neq y\wedge R(y,x,z)))$$

Next we make sure that every two vertex pairs have ranks that are comparable by containment. This means that for every $(x_1,y_1)$ and $(x_2,y_2)$ either $r_{x_1,y_1}\subseteq r_{x_2,y_2}$ or $r_{x_2,y_2}\subseteq r_{x_1,y_1}$.

$$\phi_{\mathrm{comp}}=\forall_{x_1,y_1,x_2,y_2}\neg\exists_{z_1,z_2}(R(x_1,y_1,z_1)\wedge\neg R(x_2,y_2,z_1)\wedge R(x_2,y_2,z_2)\wedge\neg R(x_1,y_1,z_2))$$

Finally, we want every non-loop vertex pair to have a non-empty rank, and moreover for it to include the constant $a$. This is crucial, because $a$ will eventually serve as an ``anchor'' making sure that there is only one way to assign ranks when encoding a full iterated matching sequence using the ternary relation $R$.

$$\phi_{\mathrm{full}}=\forall_{x,y}((x\neq y)\to R(x,y,a))$$

It is a good time to sum up the full statement that sets up our pair ranks.

$$\phi_{\mathrm{rank}}=\phi_{\mathrm{graph}}\wedge\phi_{\mathrm{comp}}\wedge\phi_{\mathrm{full}}$$

Whenever this statement is satisfied, we can use it to construct expressions that compare ranks. We will use the following expressions, which compare the ranks of $(x_1,y_1)$ and $(x_2,y_2)$,  when we formulate further conditions on $R$ that will eventually force it to conform to a full iterated matching sequence. Note that conveniently, these comparison expression also work against loops (whose ``rank'', the empty set, is considered to be the lowest).

\begin{eqnarray*}
\phi_{=}(x_1,y_1,x_2,y_2) & = & \forall_z(R(x_1,y_1,z)\leftrightarrow R(x_2,y_2,z)) \\
\phi_{\leq}(x_1,y_1,x_2,y_2) & = & \forall_z(R(x_1,y_1,z)\to R(x_2,y_2,z)) \\
\phi_{<}(x_1,y_1,x_2,y_2) & = & \phi_{\leq}(x_1,y_1,x_2,y_2)\wedge\neg\phi_{=}(x_1,y_1,x_2,y_2)
\end{eqnarray*}

\subsection{Making the ordered pairs correspond to an iterated matching}\label{ss:tm}

In this subsection we consider a ternary relation $R$ that is known to satisfy $\phi_{\mathrm{rank}}$ as defined in Subsection \ref{ss:to}, and impose further conditions that will force it to correspond to an iterated matching sequence (which will also be full by virtue of every pair having a rank).

For every rank appearing in $R$, that is for every set $A$ which is equal to $r_{x,y}$ for some $x,y\in [n]$, we refer to the set of vertex pairs having this rank as $E_i$, where $i$ is the number of ranks that appear in $R$ (including the empty set, which is the ``rank'' of loops) and are strictly contained in $A$. So in particular $E_0=\{(x,x):x\in [n]\}$, and $E_1$ for example would be the set of vertex pairs that have the smallest non-empty set as their ranks.

We first impose the restriction that for any $i$, the graph defined by $\bigcup_{j=1}^iE_j$ is a transitive graph, that is a disjoint union of cliques. By Observation \ref{ob:im} this is a necessary condition for $\bar E$ to be an iterated matching sequence (note that allowing also the $0$-ranked loops does not change the condition). This is the same as saying that for any three vertices $x,y,z$, it cannot be the case that the rank of $(x,z)$ is larger than the maximum ranks of $(x,y)$ and $(y,z)$.

$$\phi_{\mathrm{trans}}=\forall_{x,y,z}(\phi_{\leq}(x,z,x,y)\vee\phi_{\leq}(x,z,y,z))$$

Whenever $R$ satisfies the above, it is not hard to add the restriction that $E_i$ consists of disjoint complete bipartite graphs such that each of them connects exactly two components of $\bigcup_{j=1}^{i-1}E_j$, with all such components being covered. First we state that if some rank $A$ exists, that is, there exists some $(x,y)$ for which $A=r_{x,y}$, then every vertex $z$ is is part of an edge with such rank.

$$\phi_{\mathrm{cover}}=\forall_{x,y}\forall_z\exists_w\phi_{=}(x,y,z,w)$$

Then, using the prior knowledge that all connected components of both $\bigcup_{j=1}^{i-1}E_j$ and $\bigcup_{j=1}^iE_j$ are cliques, to make sure that every connected component of $E_i$ is exactly a bipartite graph encompassing two components of $\bigcup_{j=1}^{i-1}E_j$, it is enough to state that it contains no triangles, excluding of course ``triangles'' of the type $(x,x,x)$.

$$\phi_{\mathrm{part}}=\forall_{x,y,z}((x\neq y)\to\neg(\phi_{=}(x,y,y,z)\wedge\phi_{=}(x,y,x,z)))$$

All of the above is sufficient to guarantee that the relation $R$ corresponds to a full iterated matching sequence. However, as things stand now there can be many relations that correspond to the same iterated matching. This occurs because we still have unwanted freedom in choosing the sets that correspond to the possible ranks. To remove this freedom, we now require that the rank of every pair $(x,y)$ for $x\neq y$ consists of exactly one connected component of the union of the lower ranked pairs. This will be sufficient, because by $\phi_{\mathrm{full}}$ the only option for the rank would be the connected component that contains the constant $a$.

Noting that by $\phi_{\mathrm{trans}}$ these components are cliques, it is enough to require that every member of $r_{x,y}$ is connected via a lower rank edge to $a$, while every vertex that is connected to a member of $r_{x,y}$ via a lower rank edge is also a member of $r_{x,y}$. We obtain the following statement.

$$\phi_{\mathrm{anchor}}=\forall_{x,y,z}(R(x,y,z)\to(\phi_{<}(z,a,x,y)\wedge\forall_w(\phi_{<}(z,w,x,y)\to R(x,y,w))))$$

The final statement that counts the number of full iterated matching sequences, and hence provides the example proving Theorem \ref{th:tce} is the following.

$$\phi_{\mathcal M}=\phi_{\mathrm{rank}}\wedge\phi_{\mathrm{trans}}\wedge\phi_{\mathrm{cover}}\wedge\phi_{\mathrm{part}}\wedge\phi_{\mathrm{anchor}}$$

\subsection{Adapting the example to other primes}

We show here how to adapt the FOL sentence from Theorem \ref{th:tce} to provide a sequence that is not ultimately periodic modulo $p$ for any prime number $p\geq2$. The analogous corollary about removing the constant also follows.

\begin{theorem}[Ternary relation counterexample for $p\geq2$]\label{th:tcp}
For any prime number $p$, there exists an FOL sentence $\phi_{{\mathcal M}_p}$ over the language $(a,R)$, where $a$ is a (hard-coded) constant and $R$ is a relation of arity $3$, so that the corresponding class ${\mathcal C}_p$ satisfies $f_{{\mathcal C}_p}(n-1)=0$ for every $n$ that is not a power of $p$, and $f_{{\mathcal C}_p}(n-1)\equiv 1\pmod{p}$ for $n=p^m$ for every $m\in\mathbb{N}$. In particular, $f_{{\mathcal C}_p}$ is not ultimately periodic modulo $p$.
\end{theorem}

The construction follows the same lines as the extension from $p=2$ to $p\geq2$ in previous works. For completeness we give some details on how it works with respect to the version of \cite{spadd}. The basic idea is to use a ``matching'' of $p$-tuples instead of pairs.

\begin{definition}
A {\em $p$-matching} over the vertex set $[n]$ is a spanning graph, each of whose connected components is either a clique with $p$ vertices or a single vertex. A {\em perfect $p$-matching} is a $p$-matching in which there are no single vertex components (in other words, it is a partition of $[n]$ into sets of size $p$).
\end{definition}

The following is not hard to prove.

\begin{lemma}\label{lm:onepm}
There are no perfect $p$-matchings over $[n]$ unless $n$ is a multiple of $p$, in which case their number is congruent to $1$ modulo $p$.
\end{lemma}

\begin{proof}
The case where $n$ is not a multiple of $p$ is trivial. Otherwise, consider the number of possible partitions of the set $[p]$ to a sequence of subsets of sizes $i_1,\ldots,i_r$, where $\sum_{k=1}^ri_k=p$. Note that unless $i_1=p$ (and hence $r=1$), the number of such partitions is divisible by $\binom{p}{i_1}$, which is divisible by $p$ (since $p$ is a prime).

Denoting by $f_{M_p}(n)$ the number of perfect $p$-matchings over $[n]$, We consider for any $p$-matching its restriction to $[p]$ (which corresponds to a partition of $[p]$ -- the reason we need to consider the partitions as sequences rather than as unordered families of sets is that we need to consider which sets in the restriction of the $p$-matching over $[n]\setminus[p]$ they are ``attached'' to). This implies that $f_{M_p}(n)\equiv f_{M_p}(n-p)\pmod{p}$ for every $n>p$, allowing us to prove by induction that $f_{M_p}(n)\equiv 1\pmod{p}$ if $p$ divides $n$.
\end{proof}

The definition of an iterated $p$-matching sequence is what one would expect.

\begin{definition}[Iterated $p$-matching sequence]
Given a set $V$ of vertices, An {\em iterated $p$-matching sequence} is a sequence graphs over $V$, identified by their edge sets $\bar E=E_1,\ldots,E_{\ell_{\bar E}}$, satisfying the following for every $1\leq i\leq \ell_{\bar E}$.
\begin{itemize}
\item The connected components of $E_i$ are (vertex-disjoint) complete $p$-partite graphs.
\item The $p$ vertex classes of every complete $p$-partite graph in $E_i$ as above are $p$ connected components of $\bigcup_{j=1}^{i-1}E_j$ (for $i=1$ this means that $E_1$ is a $p$-matching).
\item Every connected component of $\bigcup_{j=1}^{i-1}E_j$ is a vertex class of some $p$-partite graph of $E_i$ (so in particular $E_1$ is a perfect $p$-matching).
\end{itemize}
An iterated matching sequence $\bar E$ is {\em full} if every vertex pair $u,v\in V$ (where $u\neq v$) appears in some $E_i$.
\end{definition}

Again we have the following properties, analogous to those of iterated matching sequences.

\begin{observation}
For an iterated $p$-matching $\bar E$, every $E_i$ corresponds to a perfect $p$-matching over the set of connected components of $\bigcup_{j=1}^{i-1}E_j$. Additionally, every connected component of $\bigcup_{j=1}^iE_j$ is a clique with exactly $p^i$ vertices.
\end{observation}

The above implies that there can be a full iterated matching sequence over $[n]$ if and only if $n$ is a power of $p$, in which case $\ell_{\bar E}=\log_p(n)$. Denoting the number of possible full iterated matching sequences over $[n]$ by $f_{{\mathcal M}_p}(n)$, note the following lemma.

\begin{lemma}
For every $n$ that is not a power of $p$ we have $f_{\mathcal M}(n)=0$, while for $n=p^m$ for every $m\in\mathbb N$ we have $f_{{\mathcal M}_p}(n)\equiv 1\pmod{p}$.
\end{lemma}

\begin{proof}
The case where $n$ is not a power of $p$ was already discussed above. The case $n=p^m$ is proved by induction over $m$ using Lemma \ref{lm:onepm}.
\end{proof}

From here on the construction of $\phi_{{\mathcal M}_p}$ is identical to that of $\phi_{\mathcal M}$ in Subsection \ref{ss:to} and Subsection \ref{ss:tm}, with the only exceptions being the replacements for $\phi_{\mathrm{cover}}$ and $\phi_{\mathrm{part}}$.

To construct $\phi_{\mathrm{cover}_p}$, we need to state that for every existing rank, each vertex is part of a size $p$ clique consisting of edges from this rank.

$$\phi_{\mathrm{cover}_p}=\forall_{x,y}\forall_{z_1}\exists_{z_2,\ldots,z_p}\bigwedge_{1\leq i<j\leq p}\phi_{=}(x,y,z_i,z_j)$$

To construct $\phi_{\mathrm{part}_p}$, we need to state that no $E_i$ may contain a clique with $p+1$ vertices.

$$\phi_{\mathrm{part}_p}=\forall_{z_1,\ldots,z_{p+1}}((z_1\neq z_2)\to\neg(\bigwedge_{1\leq i<j\leq p+1}\phi_{=}(z_1,z_2,z_i,z_j))$$

The final expression is the following.

$$\phi_{{\mathcal M}_p}=\phi_{\mathrm{rank}}\wedge\phi_{\mathrm{trans}}\wedge\phi_{\mathrm{cover}_p}\wedge\phi_{\mathrm{part}_p}\wedge\phi_{\mathrm{anchor}}$$

\section{A class with only one ternary relation which is not MC-finite}\label{s:s}

Starting with Theorem \ref{th:tce}, to remove the hard-coded constant (and arrive at Corollary \ref{cor:ptce}) we only need to use a single invocation of Lemma \ref{lm:hscr}. Since we started out with a single ternary relation $R$ in the language of $\phi_{\mathcal M}$, this leaves us with a statement $\phi'_{\mathcal M}$ utilizing the eight relations $R_{\emptyset},R_{\{1\}},R_{\{2\}},R_{\{3\}},R_{\{1,2\}},R_{\{2,3\}},R_{\{1,3\}},R_{\{1,2,3\}}$. In the following we show how to remove all of these relations except the relation $R_{\emptyset}$ from the language, while keeping the model counts, which leaves us with a single ternary relation. We note that the exact same treatment will also work for the modulo $p$ version $\phi'_{{\mathcal M}_p}$.

For convenience, we let $\phi_{\mathcal M,8}$ denote $\phi'_{\mathcal M}$, where ``$8$'' is the number of relations in the language. Each time we will define an expression with a smaller number of relations, and claim that the number of satisfying models is preserved.

All throughout, we assume that $\mathfrak M$ is a model for which $\mathfrak M\models\phi_{\mathcal M}$ over the language $(a,M)$ and the universe $[n]$ (where the constant $a$ is hard-coded to refer to $n$), and that $\mathfrak N$ is its corresponding model over the language of the expression $\phi_{\mathcal M,i}$ under discussion and the universe $[n-1]$ (which does not include $a$).

Referring to $\phi_{\mathrm{graph}}$, which is a component of $\phi_{\mathcal M}$, our first observation is a very easy one.

\begin{observation}
It is never the case that $\mathfrak N\models R_{\{1,2,3\}}()$, since by $\mathfrak M\models\phi_{\mathrm{graph}}$ it is never the case that $\mathfrak M\models R(a,a,a)$. Similarly, it is never the case that $\mathfrak N\models R_{\{1,2\}}(x)$ for any $x\in [n-1]$ (the universe of $\mathfrak N$), since it is never the case that that $\mathfrak M\models R(a,a,x)$.
\end{observation}

This allows us to get rid of the nullary relation $R_{\{1,2,3\}}$ and the unary relation $R_{\{1,2\}}$.

\begin{definition}
To construct $\phi_{\mathcal M,6}$ while preserving the model count, we replace all atomic formulas ``$R_{\{1,2,3\}}()$'' and ``$R_{\{1,2\}}(x)$'' (for any variable $x$) in $\phi_{\mathcal M,8}$ with the Boolean ``false'' statement $\bot$, and remove the symbols $R_{\{1,2,3\}}$ and $R_{\{1,2\}}$ from the language of $\phi_{\mathcal M,6}$.
\end{definition}

We next deal with the other two unary relations, $R_{\{1,3\}}$ and $R_{\{2,3\}}$. Here it is very important to note that the universe of $\mathfrak N$ does not include $n$, so in particular $\mathfrak N\models R_{\{1,3\}}(x)$ if and only if $\mathfrak M\models R(a,x,a)$, where $x$ is guaranteed to be unequal to $a$.

\begin{observation}
It is always the case that $\mathfrak N\models R_{\{1,3\}}(x)$ and $\mathfrak N\models R_{\{2,3\}}(x)$ for any $x$ in the universe of $\mathfrak N$, since by $\mathfrak M\models\phi_{\mathrm{full}}$ it is always the case that $\mathfrak M\models R(a,x,a)$ and $\mathfrak M\models R(x,a,a)$.
\end{observation}

This allows us to get rid of the two remaining unary relations.

\begin{definition}
To construct $\phi_{\mathcal M,4}$ while preserving the model count, we replace all atomic formulas ``$R_{\{1,3\}}(x)$'' and ``$R_{\{2,3\}}(x)$'' (for any variable $x$) in $\phi_{\mathcal M,6}$ with the Boolean ``true'' statement $\top$, and remove the symbols $R_{\{1,3\}}$ and $R_{\{2,3\}}$ from the language of $\phi_{\mathcal M,4}$.
\end{definition}

We next consider the binary relation $R_{\{3\}}$. While the truth value of $R_{\{3\}}(x,y)$ in $\mathfrak N$ depends on the actual values of $x$ and $y$, it is still fully determined by $\mathfrak N$ satisfying $\phi_{\mathcal M,4}$, or equivalently, by $\mathfrak M$ satisfying $\phi_{\mathcal M}$, because it corresponds to the truth value of $R(x,y,a)$ in $\mathfrak M$.

\begin{observation}
By $\mathfrak M\models\phi_{\mathrm{full}}$ and $\mathfrak M\models\phi_{\mathrm{graph}}$, it is always the case that $\mathfrak N\models R_{\{3\}}(x,y)$ if and only if $x\neq y$.
\end{observation}

This allows us to remove $R_{\{3\}}$ and replace it with the equivalent expression.

\begin{definition}
To construct $\phi_{\mathcal M,3}$ while preserving the model count, we replace all atomic formulas ``$R_{\{3\}}(x,y)$'' (for any variables $x$ and $y$) in $\phi_{\mathcal M,4}$ with the expression ``$(x\neq y)$'', and remove the symbol $R_{\{3\}}$ from the language of $\phi_{\mathcal M,3}$.
\end{definition}

We now consider the two remaining binary relations, $R_{\{1\}}$ and $R_{\{2\}}$. Their interpretation by $\mathfrak N$ can vary among different models satisfying $\phi_{\mathcal M,3}$. However, we note that $\mathfrak N\models R_{\{1\}}(x,y)$ if and only if $\mathfrak M\models R(a,x,y)$, and similarly $\mathfrak N\models R_{\{2\}}(x,y)$ if and only if $\mathfrak M\models R(x,a,y)$. This means that $R_{\{1\}}$ and $R_{\{2\}}$ have identical interpretations.

\begin{observation}
By $\mathfrak M\models\phi_{\mathrm{graph}}$, it is always the case that $\mathfrak N\models R_{\{1\}}(x,y)$ if and only if $\mathfrak N\models R_{\{2\}}(x,y)$ for every $x$ and $y$.
\end{observation}

This means that we can at least get rid of $R_{\{2\}}$.

\begin{definition}
To construct $\phi_{\mathcal M,2}$ while preserving the model count, we replace all atomic formulas ``$R_{\{2\}}(x,y)$'' (for any variables $x$ and $y$) in $\phi_{\mathcal M,3}$ with ``$R_{\{1\}}(x,y)$'', and remove the symbol $R_{\{2\}}$ from the language of $\phi_{\mathcal M,2}$.
\end{definition}

For the final step we cannot replace instances of $R_{\{1\}}$ with a fixed expression. However, we can ``re-purpose'' part of $R_{\emptyset}$ to hold the information currently held by $R_{\{1\}}$, allowing us to create an expression over a language containing only this one ternary relation.

For this we first note (by $\phi_{\mathrm{graph}}$) that it is never the case that $\mathfrak M\models R(x,x,y)$ for any $x,y\in [n-1]$, and hence it is never the case that $\mathfrak N\models R_{\emptyset}(x,x,y)$ whenever $\mathfrak N\models\phi_{\mathcal M,2}$. For our final transformation, we need to simulate the ``old'' $R_{\emptyset}$ using only the truth values of $R_{\emptyset}(x,y,z)$ for $x\neq y$, and then we can replace instances of $R_{\{1\}}$ by the truth values of $R_{\emptyset}(x,y,z)$ for $x=y$.

This leads us to the following definition.

\begin{definition}
To construct the final $\phi_{\mathcal M,1}$ while preserving the model count, we replace all atomic formulas ``$R_{\emptyset}(x,y,z)$'' (for any variables $x,y,z$) in $\phi_{\mathcal M,2}$ with ``$((x\neq y)\wedge R_{\emptyset}(x,y,z))$'', replace all atomic formulas ``$R_{\{1\}}(x,y)$'' in $\phi_{\mathcal M,2}$ with ``$R_{\emptyset}(x,x,y)$'', and remove the symbol $R_{\{1\}}$ from the language of $\phi_{\mathcal M,1}$.
\end{definition}

The expression $\phi_{\mathcal M,1}$, over the language containing only $R_{\emptyset}$, yields the following theorem. It is formulated modulo $2$, although as noted above it can be extended to any prime $p\geq 2$.

\begin{theorem}[A sentence with a single relation]
\label{th:main-2}
There exists an $\FOL$-sentence $\phi_{\mathcal M,1}$ over a language consisting of a single relation of arity $3$, 
so that for the class $\cC$ corresponding to $\phi_{\mathcal M,1}$, its counting function $f_{\mathcal C}(n)$ is not not ultimately periodic modulo $2$.
\end{theorem}

\part{\Large Epilogue}
\label{part-4}
\section{Conclusions and open problems}
\label{se:conclu}

In this paper we have extended the Specker-Blatter Theorem
to classes of $\tau$-structures definable in $\CMSOL$
for vocabularies $\tau$ which contain a finite number of constants, unary and binary relation symbols,
Corollary \ref{cr:spc}.
We have also shown that this does not hold already when $\tau$ consists of only one ternary relation symbol,
Theorem \ref{th:main-2}.
We also note that in \cite{pr:FischerMakowsky03,fischer2011application} we have shown that for $\cC$ definable in $\CMSOL$
such that all structures have degree bounded by a constant $d$, $S_{\cC}(n)$ is always MC-finite.
The degree of a structures $\cA$ is defined via the Gaifman graph of $\cA$.
With this the MC-finiteness of $S_{\cC}(n)$ for $\CMSOL$-definable classes of $\tau$-structures
as a function of $\tau$ is completely understood.
Applications of our results in this paper to restricted Bell numbers and various restricted partition functions
are given in \cite{ar:FischerMakowskyRakita22}.

A sequence of integers $s(n)$ is MC-finite 
if for every $m \in \N^+$ there are constants $r(m), p(m)\in \N^+$ and coefficients $\alpha_1(m), \ldots,\alpha_{p(m)} \in \N^+$,
such that for all $n \geq r(m)$ we have
$$
s(n+p(m)+1) \equiv  \sum_{i=0}^{p(m)} \alpha_i s(i) \mod{m}
$$
The Specker-Blatter Theorem gives little information on the constants 
$r(m), p(m)$ or the coefficients $\alpha_1(m), \ldots, \alpha_{p(m)}$.
These in particular depend on the substitution rank of the class $\cC$.
In fact Theorem \ref{th:rank} gives a very bad estimate of the substitution rank in the case of binary relation symbols.
The constants are computable, but it is not known whether they are always computable in feasible time
or whether their size is bounded by an elementary function.
In the presence of constants the substitution rank is not defined. Our main Theorem \ref{th:r} allows to eliminate
the constants, and therefore gives a formula for which the substitution rank is defined. However,
due to the increased complexity of the resulting formula, the estimate of the substitution rank
will be even worse.

\begin{problem}
Given a sentence $\phi$ in $\CMSOL(\tau)$ where $\tau$ consists only of constants, unary and binary relation symbols,
\begin{enumerate}[(i)]
\item
what is the time complexity of computing the constants $r(m), p(m)$ and the coefficients $\alpha_1(m), \ldots, \alpha_{p(m)}$?
\item
what can we say about the size of these constants?
\end{enumerate}
\end{problem}

The proof of Theorem \ref{th:rank} depends on the Feferman-Vaught Theorem which also holds for
$\CMSOL(\tau)$ for any finite relational $\tau$, \cite{ar:FV,ar:MakowskyTARSKI}.
In our context, the Feferman-Vaught Theorem allows to check whether a formula of $\CMSOL(\tau)$
holds in $Subst(\cA_1,a,\cA_2)$ by checking a sequence of $\CMSOL(\tau)$-formulas in $\cA_1$ and $\cA_2$ 
independently. This sequence is called a reduction sequence, cf. \cite{sbplus}.
In \cite{dawar2007model} it is shown that even for $\FOL(\tau)$ the size of the reduction sequences
for the Feferman-Vaught Theorem cannot, in general, be bounded by an elementary functions.

\begin{problem}
Does there exist an elementary function $f(k)$, so that for any sentence $\phi$ in $\CMSOL(\tau)$ 
where $\tau$ consists only of constants, unary and binary relation symbols, 
the size of the constants $r(m)$ and $p(m)$ is bounded by $f(\max\{|\phi|,m\})$?
\end{problem}


The Specker-Blatter Theorem also applies to hereditary, monotone  and minor-closed graph classes,
provided they are definable using a finite set of forbidden (induced) subgraphs or minors.
In the first two cases such a class is $\FOL$-definable. In the case of a minor-closed class,
B. Courcelle showed that it is $\MSOL$-definable, see \cite{bk:CourcelleEngelfriet2011}.
By the celebrated theorem of N. Robertson and P. Seymour, \cite{bk:Diestel05}, every minor-closed class of graphs
is definable by a finite set of forbidden minors. However, there are monotone (hereditary)
classes of graphs where a finite set of forbidden (induced) subgraphs does not suffice.

\begin{problem}
Are there hereditary  or monotone classes of graphs $\cC$ such that $Sp_{\cC}(n)$ is not MC-finite?
\end{problem}

An analogue question arises when we replace graphs by finite relational $\tau$-structures. 
In this case one speaks of classes of $\tau$-structures closed under substructures.
Every class of finite $\tau$-structures $\cC$ closed under substructures can be characterized by a set of
forbidden substructures. If this set is finite, $\cC$ is again $\FOL$-definable, and the Specker-Blatter
Theorem applies.

\begin{problem}
\begin{enumerate}[(i)]
\item
Let $\tau$ be a relational vocabulary. 
Are there  substructure closed classes $\cC$ of $\tau$-structures such that $Sp_{\cC}(n)$ is not MC-finite?
\item
Same question when all the
relations are at most binary? 
\end{enumerate}
\end{problem}

\subsection*{Acknowledgments}
We would like to thank J. Baldwin for useful comments.

\addcontentsline{toc}{section}{References}
\bibliographystyle{plain}
\bibliography{paper-ref}
\end{document}